\documentclass[a4paper,12pt]{amsart}
\usepackage{amssymb}
\usepackage{amsmath}

\def\.{\cdot}

\def\a{\alpha}
\def\b{\beta}
\def\c{\gamma}

\def\vs{\vskip .6cm}

\def\n{\nabla}

\def\l{\lambda}
\def\s{\sigma}
\def\t{\tilde}
\def\beq{\begin{equation}}
\def\eeq{\end{equation}}
\def\bea{\begin{eqnarray*}}
\def\eea{\end{eqnarray*}}
\def\beaa{\begin{eqnarray}}
\def\eeaa{\end{eqnarray}}
\def\ba{\begin{array}}
\def\ea{\end{array}}
\def\f{\varphi}
\def\ff{\phi}
\def\o{\omega}

\def\L{\Lambda}
\def\G{\Gamma}
\def\bp{\begin{proof}}
\def\r{\end{proof}}
\def\res{\arrowvert}
\def\V{{\mathcal V}}
\def\A{{\mathcal A}}

\def \RM{\mathbb{R}}
\def \Ca{\mathbb{O}}

\def \ZM{\mathbb{Z}}
\def \CM{\mathbb{C}}

\def \HM{\mathbb{H}}
\def \SM{\mathbb{S}}
\def \PM{\mathbb{P}}

\def\d{{\delta}}

\def\del{{\partial}}
\def\Ric{\mathrm{Ric}}
\def\id{\mathrm{id}}
\def\be{\begin{equation}}
\def\ee{\end{equation}}
\def\tr{\mathrm{tr}}

\def\Hol{\mathrm{Hol}}

\def\hol{\mathfrak{hol}}
\def\so{\mathfrak{so}}
\def\su{\mathfrak{su}}
\def\spin{\mathfrak{spin}}
\def\sp{\mathfrak{sp}}
\def\gg{\mathfrak{g}}
\def\hh{\mathfrak{h}}
\def\pp{\mathfrak{p}}

\def\k{\mathfrak{k}}
\def\kk{\kappa}

\def\pp{\mathfrak{p}}

\def\SU{\mathrm{SU}}
\def\U{\mathrm{U}}
\def\E{\mathrm{E}}
\def\F{\mathrm{F}}
\def\SO{\mathrm{SO}}

\def\End{\mathrm{End}}
\def\Iso{\mathrm{Iso}}
\def\PSO{\mathrm{PSO}}
\def\Cl{\mathrm{Cl}}
\def\Sp{\mathrm{Sp}}
\def\Spin{\mathrm{Spin}}
\def\Pin{\mathrm{Pin}}
\def\Ker{\mathrm{Ker}}
\def\ad{\mathrm{ad}}
\def\res{\arrowvert}

\def\rr{r}

\def\scal{\mathrm{scal}}
\def\Ad{\mathrm{Ad}}

\def\I{{\mathcal I}}

\newtheorem{epr}{Proposition}[section]
\newtheorem{ath}[epr]{Theorem}
\newtheorem{elem}[epr]{Lemma}
\newtheorem{ecor}[epr]{Corollary}

\theoremstyle{definition}

\newtheorem{ede}[epr]{Definition}
\newtheorem{ere}[epr]{Remark}
\newtheorem{exe}[epr]{Example}

\title{Clifford structures on Riemannian manifolds}

\author{Andrei Moroianu and Uwe Semmelmann}

\address{Andrei Moroianu \\ CMLS\\ {\'E}cole Polytechnique \\ UMR 7640 du CNRS
\\ 91128 Palaiseau \\ France}
\email{am@math.polytechnique.fr}

\address{Uwe Semmelmann\\ 
Institut f\"ur Geometrie und Topologie \\
Fachbereich Mathematik\\
Universit{\"a}t Stuttgart\\
Pfaffenwaldring 57 \\
70569 Stuttgart, Germany}
\email{uwe.semmelmann@mathematik.uni-stuttgart.de}

\thanks{This work was supported by the French-German cooperation
  project Procope no. 17825PG}

\begin{document}

\begin{abstract}
We introduce the notion of {\em even Clifford structures} on
Riemannian manifolds, which for rank $\rr=2$ and $\rr=3$ reduce to
almost Hermitian and quaternion-Hermitian structures respectively. We give the
complete classification of manifolds carrying {\em parallel} rank $\rr$
even Clifford structures: K\"ahler, quaternion-K\"ahler and Riemannian
products of 
quaternion-K\"ahler manifolds for $\rr=2,\,3$ and $4$ respectively,
several classes of $8$-dimensional manifolds (for
$5\le\rr\le 8$), families of real, complex and
quaternionic Grassmannians (for $\rr=8,\ 6$ and $5$
respectively), and Rosenfeld's elliptic projective planes $\Ca \PM^2$,
$(\CM\otimes
\Ca) \PM^2$, $(\HM\otimes \Ca) \PM^2$ and $(\Ca\otimes \Ca) \PM^2$, which
are symmetric spaces associated to
the exceptional simple Lie groups $\F_4$, $\E_6$, $\E_7$ and $\E_8$ (for
$\rr=9,\ 10,\ 12$ and $16$ respectively). As an application, we
classify all Riemannian manifolds whose metric is bundle-like along the
curvature constancy distribution,
generalizing well known results in Sasakian and 3-Sasakian geometry.

\bigskip

\noindent
2000 {\it Mathematics Subject Classification}: Primary 53C26,
53C35, 53C10, 53C15.

\medskip
\noindent{\it Keywords:} Clifford structure, K\"ahler,
quaternion-K\"ahler, symmetric spaces, exceptional Lie groups,
Rosenfeld's elliptic projective planes,
curvature constancy, fat bundles.
\end{abstract}

\maketitle

\section{Introduction}

The main goal of the present paper is to introduce a new
algebraic structure on Riemannian manifolds, which we refer to as 
{\em Clifford structure}, 
containing almost complex structures and almost quaternionic structures
as special cases. 

Roughly speaking, by a
Clifford (resp. even Clifford) structure on a Riemannian manifold $(M,g)$ 
we understand a Euclidean vector 
bundle $(E,h)$ over $M$, called {\em Clifford bundle}, 
together with a representation of the 
Clifford algebra bundle $\Cl(E,h)$ (resp. $\Cl^0(E,h)$) on 
the tangent bundle $TM$. One might notice the duality between 
spin and Clifford structures: While in spin geometry, the spinor bundle is 
a representation space 
of the Clifford algebra bundle of $TM$, in the new framework, it is the
tangent bundle of the manifold which becomes a representation space of the
(even) Clifford algebra bundle of the Clifford bundle $E$.

Several approaches to the concept of Clifford structures on
Riemannian manifolds can be found in the literature. We must stress 
from the very beginning
on the somewhat misleading fact that the same terminology 
is used for quite different notions.
Most authors have introduced Clifford structures 
as a family of {\em global} almost complex 
structures satisfying the Clifford relations, i.e. as a pointwise
representation of the Clifford algebra $\Cl_n$ on each tangent space
of the manifold. In the sequel we will refer to these structures as
{\em flat Clifford 
structures}. In contrast, our definition only involves
{\em local} almost complex structures, obtained from local orthonormal
frames of the Clifford bundle $E$, and reduces to the previous notion
when $E$ is trivial.

Flat Clifford structures were considered by Spindel et al. in
\cite{spindel}, motivated by the fact that in the 2-dimensional
supersymmetric $\sigma$-model, a target 
manifold with $N-1$ independent parallel anti-commuting complex
structures gives rise to $N$ supersymmetries. They claimed that 
on group manifolds $N \le 4$ but later on, Joyce showed that this
restriction does not hold in the
non-compact case (cf. \cite{joyce}) and provided a method 
to construct manifolds with arbitrarily large Clifford
structures. At the same time, Barberis et al. constructed in
\cite{barberis} flat Clifford structures on compact flat manifolds, 
by means of 2-step nilpotent Lie groups. 

Yet another notion of Clifford
structures was used in connection with the Osserman Conjecture.
Following ideas of Gilkey, Nikolayevsky  defined in
\cite{niko} Clifford structures on Riemannian manifolds with an
additional assumption on the Riemannian curvature tensor. 

An author who comes close to our concept of even Clifford structure, 
but restricted to a particular case, is Burdujan. His 
{\em Clifford-K\"ahler} manifolds, introduced 
in \cite{burdujan1} and \cite{burdujan2}, 
correspond in our terminology to manifolds 
with a rank 5 parallel even Clifford structure. He
proves that such manifolds have to be Einstein (a special case of
Proposition \ref{p37} below). Note also that Spin(9)-structures on
16-dimensional manifolds studied by Friedrich \cite{fr} correspond to rank
9 even Clifford structures in our setting.

The core of the paper consists of the classification of manifolds carrying 
parallel even Clifford structures, cf. Theorem \ref{class}. 
In rank $\rr=2$ and $\rr=3$ this reduces to 
K\"ahler and quaternion-K\"ahler structures respectively.
We obtain Riemannian products of
quaternion-K\"ahler manifolds for $\rr=4$,
several classes of $8$-dimensional manifolds (for
$5\le\rr\le 8$), families of real, complex and
quaternionic Grassmannians (for $\rr=8,\ 6$ and $5$
respectively), and Rosenfeld's elliptic projective planes $\Ca \PM^2$,
$(\CM\otimes
\Ca) \PM^2$, $(\HM\otimes \Ca) \PM^2$ and $(\Ca\otimes \Ca) \PM^2$, which
are symmetric spaces associated to
the exceptional simple Lie groups $\F_4$, $\E_6$, $\E_7$ and $\E_8$ (for
$\rr=9,\ 10,\ 12$ and $16$ respectively). Using similar arguments we
also classify  
manifolds carrying parallel Clifford structures, showing 
that parallel Clifford structures can only exist in low
rank ($\rr\le 3$), in low dimensions ($n\le 8$) or on flat
spaces (cf. Theorem \ref{t23}).

In Section \ref{scc}, we give a geometric application of our
classification theorem to the theory of manifolds with {\em curvature
  constancy}, a notion
introduced in the 60's by Gray \cite{gray}. Roughly speaking, a tangent vector
$X$ on a Riemannian manifold $(Z,g_Z)$ belongs to the curvature
constancy $\V$ if its contraction with the 
Riemannian curvature tensor $R^Z$ equals its contraction with the
algebraic curvature tensor of the round sphere, cf. \eqref{cn}
below. One reason why Gray was interested in this notion is that on
the open set of $Z$ where the dimension of the curvature constancy
achieves its minimum, $\V$ is a totally geodesic distribution whose
integral leaves are locally isomorphic to the round sphere.

Typical examples of manifolds with non-trivial curvature constancy are
Sasakian and 3-Sa\-sakian manifolds, the dimension of $\V$ being
(generically) 1 and 3 respectively. Rather curiously, Gray seems to have
overlooked these examples when he conjectured in \cite{gray} that if
the curvature constancy of a Riemannian manifold $(Z,g_Z)$ is
non-trivial, then the manifold is locally isometric to the round sphere.
By the above, this conjecture is clearly false, but one may wonder whether
counter-examples, other than Sasakian and 3-Sasakian structures, do exist.

Using Theorem \ref{class}, we classify Riemannian 
manifolds $Z$ admitting non-trivial
curvature constancy $\V$ under the additional assumption that the metric
is {\em bundle-like} along the distribution $\V$, i.e. such that
$Z$ is locally the total space of a Riemannian submersion $Z\to M$ whose
fibres are the integral leaves of $\V$, cf. Theorem \ref{fc}. Notice
that Sasakian and 3-Sasakian manifolds appear in this classification,
being total  spaces of (locally defined) Riemannian submersion over
K\"ahler and quaternion-K\"ahler manifolds respectively. 

Bundle-like metrics with curvature constancy also occur as a special
case of {\em fat bundles}, introduced by
Weinstein in \cite{weinstein} and revisited by Ziller (cf.
\cite{ziller}, \cite{z2}). A Riemannian submersion is called fat if the
sectional curvature is positive on planes spanned by a horizontal
and a vertical vector. Homogeneous fat bundles were classified
by B\'erard-Bergery in \cite{berard}. It turns out that all our homogeneous
examples with curvature constancy (cf. Table 3) may be
found in his list. It is still an open question whether 
examples of non-homogeneous fat bundles with fibre dimension larger
than one exist (cf. \cite{ziller}).

{\em Acknowledgments.} We are indebted to P. Gauduchon for several
suggestions and remarks, to P.-A. Nagy who 
has drawn our attention on the curvature constancy problem, and to W. Ziller 
for having communicated us his work \cite{ziller} on fat bundles and
for several useful comments. We also thank the anonymous referee for
his suggestions which helped us to improve the exposition.

\vs

\section{Clifford structures}

We refer to \cite{LM} for backgrounds on Clifford algebras and
Clifford bundles.

\begin{ede}\label{dcs}
A rank $\rr$ {\em Clifford structure} on a Riemannian
manifold $(M^n,g)$ is an oriented rank $\rr$ Euclidean bundle
$(E,h)$ over $M$ together with a non-vanishing algebra bundle
morphism, called {\em Clifford morphism}, $\f:\Cl(E,h)\to \End(TM)$
which maps $E$ into the bundle of 
skew-symmetric endomorphisms $\End^-(TM)$.
\end{ede}

The image by $\f$ of every unit vector $e\in E_x$ is a Hermitian
structure $J_e$ on $T_xM$ (i.e. a complex structure compatible 
with the metric $g$):
$$J_e ^2=\f(e)\circ\f(e)=\f(e\. e)=\f(-h(e,e))=-\id_{T_xM}.$$
Since the square norm of a Hermitian structure $J$ is equal to the dimension
$n$ of the space on which it acts, we see that $(E,h)$ can be identified by
$\f$ with its image $\f(E)\subset \End^-(TM)$ endowed with the 
Euclidean metric $\frac1n g$.

The universality property of the Clifford algebra immediately shows
that a rank $\rr$ Clifford structure on $(M,g)$ is a rank $\rr$
sub-bundle of $\End^-(TM)$, locally spanned by anti-commuting almost complex
structures $J_i$, $i=1,\ldots,\rr$. 

In terms of $G$-structures, a Clifford structure is equivalent to a
reduction of the orthonormal frame bundle of $M$. More precisely,
the restriction of the Clifford map $\f$ to some fibre $\Cl(E,h)_x$
defines, up to conjugacy, a representation of the Clifford algebra
$\Cl_{\rr}$ on $\RM^n$. This representation maps the groups
$\Pin(\rr)$ and $\Spin(\rr)$ isomorphically onto subgroups of
$\SO(n)$. If $C(Pin(\rr))$ denotes the centralizer of $\Pin(\rr)$ in
$\SO(n)$, then a Clifford structure is equivalent to a
reduction of the structure group of $M$ to
$\Spin(\rr)\cdot C(Pin(\rr))\subset \SO(n)$. We skip the (rather
straightforward) proof of this fact, since we will not need it in the sequel.

A Clifford structure $(M,g,E,h)$ is called {\em parallel} if the
sub-bundle $\f(E)$ of $\End^-(TM)$ is parallel with respect to the
Levi-Civita connection $\n^g$ of $(M,g)$. 

Since every oriented rank 1 vector bundle is trivial, 
there is a one-to-one correspondence between rank 1 Clifford structures and
almost Hermitian structures on $(M,g)$. 
A rank 1 Clifford structure is parallel if and only if the
corresponding almost Hermitian structure is K\"ahler.

Every hyper-K\"ahler manifold $(M^n,g,I,J,K)$ carries parallel rank 2
Clifford structures (e.g. the sub-bundle of $\End^-(TM)$ generated
by $I$ and $J$). The converse holds for $n>4$ (cf. Theorem \ref{t23}
below). Notice also that by the very definition, 
a quaternion-K\"ahler structure is nothing else but a parallel rank 3 Clifford
structure. 

The classification of $n$-dimensional Riemannian
manifolds carrying rank $\rr$
parallel Clifford structures will be given in Theorem \ref{t23} below.
It turns out that parallel Clifford structures can only exist either in low
ranks ($\rr\le 3$), or in low dimensions ($n\le 8$) or on flat
spaces. Therefore, even though it provides a common
framework for K\"ahler, quaternion-K\"ahler and hyper-K\"ahler
geometries, the notion of parallel Clifford structure is in some
sense too restrictive. 

We will now introduce a natural extension of Definition \ref{dcs}, by 
requiring the Clifford morphism to be defined only on
the {\em even} Clifford algebra bundle of $E$. We obtain in this way 
much more flexibility and examples, while a complete classification 
in the parallel case is still possible.

\begin{ede}\label{ecs}
A rank $\rr$ {\em even Clifford structure} ($\rr\ge 2$) on a Riemannian
manifold $(M^n,g)$ is an oriented rank $\rr$ Euclidean bundle
$(E,h)$ over $M$ together with an algebra bundle
morphism, called {\em Clifford morphism}, $\f:\Cl^0(E,h)\to \End(TM)$
which maps $\L^2E$ into the bundle of skew-symmetric endomorphisms
$\End^-(TM)$. Recall that $\L^2E$ is viewed as a
sub-bundle of $\Cl^0(E,h)$ by identifying $e\wedge f$ with
$e\. f+h(e,f)$ for every $e,f\in E$. Two even Clifford structures
$(E_1,h_1,\f_1)$ and $(E_2,h_2,\f_2)$ are isomorphic if there exists an
algebra bundle isomorphism $\lambda:\Cl^0(E_1,h_1)\to\Cl^0(E_2,h_2)$
such that $\f_2\circ\lambda=\f_1$.
\end{ede}

\begin{ere}\label{e}
Since the definition above only involves the exterior power $\Lambda
^2E$, the bundle $E$ itself is not part of an even Clifford
structure. As a matter of fact, there exist isomorphic even Clifford
structures with non-isomorphic bundles $E$ (see Example \ref{4.1} below).
\end{ere}

As before, an even Clifford structure is equivalent to the reduction
of the orthonormal frame bundle of $M$ to the subgroup $S\cdot C(S)$ of
$\SO(n)$, where $S$ denotes the image of
$\Spin(\rr)$ in $\SO(n)$ through the representation of
the even Clifford algebra $\Cl^0_{\rr}$ on $\RM^n$ defined (up to
conjugacy) by the map $\f$, and $C(S)$ is the centralizer of $S$ in $\SO(n)$.  
In more familiar terms, an even Clifford structure can be
characterized as follows:

\begin{elem} \label{or} Let $(E,h)$ be a rank $\rr$ even Clifford
  structure and let $\{e_1,\ldots,e_\rr\}$ be a local
  $h$-orthonormal frame on $E$.
The local endomorphisms $J_{ij}:=\f(e_i\.e_j)\in\End(TM)$
are skew-symmetric for $i\ne j$ and satisfy
\beq\label{cstr}
\begin{cases}J_{ii}=-\id\qquad& \; \mbox{for all}\; \ 1\le i\le\rr,\\
J_{ij}=-J_{ji} \mbox{ and } J_{ij}^2=-\id & \; \mbox{for all}\;\ i\ne j,\\
J_{ij}\circ J_{ik}=J_{jk}& \; \mbox{for all $i,j,k$
  mutually distinct},\\
J_{ij}\circ J_{kl}=J_{kl}\circ J_{ij}& \; \mbox{for all $i,j,k,l$
  mutually distinct}.
\end{cases}
\eeq
Moreover, if $\rr\ne 4$, then
\beq\label{ort}
<J_{ij},J_{kl}>=0,\qquad \mbox{unless } i=j,\ k=l\ \mbox{or}\ i=k\ne
j=l\ \mbox{or}\ i=l\ne k=j.
\eeq

\end{elem}

\bp 
The first statements follow directly from the usual relations in
the Clifford algebra
$$\begin{cases}e_i\. e_j=-e_j\.e_i\qquad& \; \mbox{for all}\; \ i\ne j,\\
(e_i\.e_j)^2=-\id & \; \mbox{for all}\;\ i\ne j,\\
(e_i\.e_j)\.(e_i\.e_k)=(e_j\.e_k)& \; \mbox{for all $i,j,k$
  mutually distinct},\\
(e_i\.e_j)\.(e_k\.e_l)=(e_k\.e_l)\.(e_i\.e_j)& \; \mbox{for all $i,j,k,l$
  mutually distinct}.
\end{cases}$$
The orthogonality of $J_{ij}$ and $J_{kl}$ is obvious when exactly two
of the subscripts coincide (since the corresponding endomorphisms
anti-commute). For $\rr=3$, \eqref{ort} is thus satisfied. Assume now
that $\rr\ge 5$ and that all four subscripts are mutually distinct. We
then choose $s$ different from $i,j,k,l$ and write, using the fact
that $J_{sl}$ and
$J_{ij}$ commute:
\bea<J_{ij},J_{kl}>&=&\tr(J_{ij}J_{kl})=\tr(J_{ij}J_{sk}J_{sl})=
\tr(J_{sl}J_{ij}J_{sk})\\
&=&\tr(J_{ij}J_{sl}J_{sk})=\tr(J_{ij}J_{lk})=-<J_{ij},J_{kl}>.
\eea
\r

Every Clifford structure $E$ induces an even Clifford
structure of the same rank. To see this, one needs to check on a local
  orthonormal frame $\{e_1,\ldots,e_\rr\}$ of $E$ that 
$\f(e_i\wedge e_j)$ is skew-symmetric for all $i\ne j$. This is due to
the fact that $\f(e_i\wedge e_j)=\f(e_i)\circ\f(e_j)$ is the composition of
two anti-commuting skew-symmetric endomorphisms.

The converse also holds if the rank of the Clifford bundle $E$ is equal to
$3$ modulo 4. Indeed, if $\rr=4k+3$, the Hodge
isomorphism $E\simeq \L^{\rr-1}E\subset \Cl^0(E,h)$ extends by the
universality property of the Clifford algebra to an
algebra bundle morphism $h:\Cl(E,h)\to\Cl^0(E,h)$. Thus, if
$\f:\Cl^0(E,h)\to\End(TM)$ is the Clifford morphism defining the even
Clifford structure, then $\f\circ h:\Cl(E,h)\to\End(TM)$ is an algebra
bundle morphism 
mapping $E$ into $\End^-(TM)$ (because the image by $\f\circ h$ of every
element of $E$ is a composition of $2k+1$ mutually commuting
skew-symmetric endomorphisms of $TM$). 

If the rank of the Clifford bundle $E$ is not equal to
$3$ modulo 4, the representation of $\Cl^0(E,h)$ on $TM$ cannot be 
extended in general to a representation of the whole Clifford algebra bundle
$\Cl(E,h)$. This can be seen on examples as follows. If $\rr=1,\,2,\,4$ or $8$
modulo $8$, one can take $M=\RM^n$ to be the representation space of an
irreducible representation of $\Cl^0_{\rr}$ and $E$ to be the trivial vector
bundle of rank $\rr$ over $M$. Then the obvious even Clifford
structure $E$ does not extend to a Clifford structure simply for
dimensional reasons (the dimension of any irreducible representation of
$\Cl_{\rr}$ is twice the dimension of any irreducible representation of
$\Cl^0_{\rr}$ for $\rr$ as above). For $\rr=5$, an example is provided
by the quaternionic projective space $\HM\PM^2$ which carries an even
Clifford structure of rank 5 (cf. Theorem \ref{class}). On the other
hand, any Riemannian manifold carrying a rank 5 Clifford structure is
almost Hermitian (with respect to the endomorphism induced by the
volume element of the Clifford algebra bundle), and it is well known that
$\HM\PM^2$ carries no almost complex structure \cite{massey} (cf. also
\cite{gms}). 
Finally, for $\rr=6$, an example is given by the
complex projective space $\CM\PM^4$, which carries a rank 6 even Clifford
structure (cf. Theorem \ref{class}), but no rank 6 Clifford structure,
since this would imply the triviality of its canonical bundle. Similar
examples can be constructed for all $\rr=5$ and $6\mod 8$.

An even Clifford structure $(M,g,E,h)$ is called {\em
  parallel}, if there exists a metric connection $\n^E$ on $(E,h)$ such
that $\f$ is connection preserving, i.e.
\beq\label{para}\f(\n^E_X\s)=\n^g_X\f(\s)\eeq for every tangent
vector $X\in TM$ and section $\s$ of $\Cl^0(E,h)$.

\begin{ere}\label{re} For $\rr$ even, the notion of an even Clifford
  structure of rank $\rr$ admits a slight
  extension to the case where $E$ is no longer a vector bundle but a
  {\em projective} bundle, i.e. a locally defined vector bundle
  associated to some $G$-principal bundle via a projective
  representation $\rho:G\to\PSO(\rr)=\SO(\rr)/\{\pm I_\rr\}$. Since
  the extension of the
  standard representation of $\SO(\rr)$ from $\RM^\rr$ to
  $\L^2\RM^\rr$ factors through $\PSO(\rr)$, we see that the second
  exterior power of any projective vector bundle is a well-defined
  vector bundle, so Definition \ref{ecs} can be adapted to this
  setting and the corresponding structure will be referred to as 
  {\em projective even Clifford structure} in the sequel.
\end{ere}

The main goal of this section is to classify (cf. Theorem \ref{class})
complete simply
connected Riemannian manifolds $(M,g)$ which carry parallel even
Clifford structures as introduced in Definition \ref{ecs}, in
the extended sense of Remark \ref{re}. The results are 
listed in Tables 1 and 2 below. The classification of manifolds
carrying parallel Clifford structures will then be obtained as
a by-product of Theorem \ref{class} by a case-by-case analysis.

We start by examining even Clifford structures of low rank.

\begin{exe}\label{4.1}
A rank $2$ even Clifford structure induces an almost Hermitian structure on
$(M,g)$ (the image by $\f$ of the volume element of
$\L^2E$). Conversely, every almost Hermitian structure $J$ on $(M,g)$
induces a rank 2 even Clifford structure by taking $(E,h)$ to be an
arbitrary oriented rank 2 Euclidean bundle (see Remark \ref{e}) and
defining $\f$ by the fact 
that it maps the volume element of $(E,h)$ onto $J$. An even Clifford
structure is parallel if and only if the corresponding almost
Hermitian structure $J$ is a K\"ahler structure on $(M,g)$.
\end{exe}

\begin{exe}\label{4.2}
A rank $3$ even Clifford structure induces a quaternionic structure on
$(M,g)$ i.e. a rank $3$ sub-bundle $S$ of $\End(TM)$ locally spanned
by three almost Hermitian structures satisfying the quaternion
relations. If $\{e_1,e_2,e_3\}$ is a local orthonormal basis of $E$,
$S$ is spanned by $I:=\f(e_1\.e_2)$, $J:=\f(e_2\.e_3)$ and
$K:=\f(e_3\.e_1)$. Conversely, every quaternionic structure $S$ on
$(M,g)$ induces a rank 3 even Clifford structure by taking $E=S$ with the
induced Euclidean structure and defining $\f$ as the Hodge isomorphism
$\L^2E\simeq E=S$. By this correspondence, a parallel even Clifford
structure is equivalent to a quaternion-K\"ahler structure on $M$.

Note that the quaternion-K\"ahler condition is empty in dimension 4.
There are several ways to see this, e.g. by saying that
$\Sp(1)\.\Sp(1)=\SO(4)$ so there is no holonomy restriction. In our
setting, this corresponds to the fact that the bundle $E:=\L^2_+M$
of self-dual 2-forms canonically defines a rank 3 parallel even
Clifford structure on every 4-dimensional (oriented) Riemannian
manifold.
\end{exe}

We thus see that K\"ahler and quaternion-K\"ahler geometries fit
naturally in the more general framework of parallel even Clifford
structures.

The isomorphism $\so(4)\simeq \so(3)\oplus\so(3)$ reduces the case
$\rr=4$ to $\rr=3$ (see Proposition \ref{p37} (i) below).

Let us now make the following:

\begin{ede} A parallel even Clifford structure $(M,E,\n^E)$ is called {\em
    flat} if the connection $\n^E$ is flat.
\end{ede}

\begin{ath} \label{t35}
A complete simply connected Riemannian manifold $(M^n,g)$
carrying a flat even Clifford
  structure $E$ of rank $\rr\ge5$ is flat, (and thus isometric with a
  $\Cl^0_\rr$ representation space).
\end{ath}

\bp One can choose a parallel global
  orthonormal frame $\{e_i\}, i=1,\ldots, \rr,$ on $E$, which induces
  global parallel 
  complex structures $J_{ij}:=\f(e_i\. e_j)$ on $M$ for every $i<j$.

  We claim that if $M$ is irreducible, then it is flat. Since $M$ is
  hyper-K\"ahler
  with respect to the triple $J_{12}$, $J_{31}$, $J_{23}$, it has to
  be Ricci-flat. According to the Berger-Simons Holonomy
  Theorem (cf. \cite{besse}, p. 300), $M$ is either symmetric (hence
  flat, since a 
  symmetric Ricci-flat manifold is flat), or has holonomy
  $\SU(n/2),\, \Sp(n/4)$ or $\Spin(7)$. The last three cases actually
  do not occur. Indeed, the
  space of parallel 2-forms on $M$ corresponds to the fixed points of
  the holonomy representation on $\L^2\RM^n$, or equivalently to the
  centralizer of the holonomy Lie algebra $\hol(M)$ in $\so(n)$.
  This centralizer is zero for $\Hol(M)=\Spin(7)$, 1-dimensional for
  $\Hol(M)=\SU(n/2)$ and 
  $3$-dimensional
  $\Hol(M)=\Sp(n/4)$. On the other hand, the space of parallel 2-forms
  on $M$ has
  dimension at least $\rr-1\ge 4$ ( any two of $J_{1i}$, $1<i\le\rr$
  anti-commute so they are linearly independent), a contradiction
  which proves our claim.

Back to the general case, the de Rham decomposition theorem states
that $M$ is a Riemannian product $M=M_0\times M_1\times\ldots\times
M_k$, where
$M_0$ is flat, and each $M_i$, $i\ge 1$ is irreducible, non-flat.
It is well known that a parallel complex structure $J$ on a
Riemannian manifold $(M,g)$ preserves the tangent bundle of every
irreducible non-flat factor of $M$. Indeed, if $M_1$ is such a factor,
then $J(TM_1)\cap TM_1$ is a parallel sub-bundle of $TM_1$, so either
$J(TM_1)=TM_1$ or $g(JX,Y)=0$ for all $X,Y\in TM_1$. But the latter
case is impossible since otherwise the Bianchi identity would imply
\bea
R(X,Y,X,Y)&=&R(X,Y,JX,JY)=R(X,JX,Y,JY)+R(JX,Y,X,JY)=0
\eea
for all  $X,Y\in TM_1$, so $M_1$ would be flat.

Consequently, each non-flat irreducible factor in the de Rham
decomposition of $M$
is preserved by every $J_{ij}$, and thus inherits a flat even Clifford
  structure of rank $\rr$. The first part of the proof shows that no
  such factor exists,
  so $M=M_0$ is flat.
\r

The next result is crucial for the classification of parallel even
Clifford structures. 

\begin{epr}\label{p37}
Assume that the complete simply connected Riemannian manifold $(M^n,g)$ carries
a parallel
non-flat even Clifford structure $(E,\n^E)$ of rank $\rr\ge 3$. Then
the following holds:
\begin{enumerate}
\item[(i)] If $\rr=4$ then $(M,g)$ is a Riemannian product of two
  quaternion-K\"ahler
manifolds.
\item[(ii)] If $\rr\ne 4$ and $n\ne 8$ then
\begin{enumerate}
\item[(a)] The curvature of $\n^E$, viewed as a map from $\L^2M$ to
  $\End^-(E)\simeq
\L^2E$ is a non-zero constant times the metric adjoint of the Clifford
map $\f$.
\item[(b)] $M$ is Einstein with non-vanishing
scalar curvature and has irreducible holonomy.
\end{enumerate}
\item[(iii)] If $\rr\ne 4$ and $n=8$, then (a) implies (b).
\end{enumerate}
\end{epr}

\bp

Any local orthonormal frame $\{e_1,\ldots,e_\rr\}$ on $E$
induces local endomorphisms on $M$ defined as before by $J_{ij}:=\f(e_i\.e_j)$.
We denote by $\o_{ij}$ the curvature forms of the connection $\n^E$
with respect to the local frame $\{e_i\}$:
$$R^E_{X,Y}e_i=\sum_{j=1}^\rr\o_{ji}(X,Y)e_j.$$
From \eqref{para} we immediately get $\;\f\circ
R^E_{X,Y}=R_{X,Y}\circ\f$, where $R$ denotes the Riemannian
curvature tensor on $(M,g)$. Consequently,
\begin{equation}\label{0r}
\begin{split}
R_{X,Y}J_{ij}=&R_{X,Y}\f(e_i\.e_j)=\f[R^E_{X,Y}(e_i\.e_j)]\\
=&\f[\sum_{s=1}^\rr\o_{si}(X,Y)e_s\.e_j+e_i\.\sum_{s=1}^\rr\o_{sj}(X,Y)e_s]\\
=&\sum_{s=1}^\rr[\o_{si}(X,Y)J_{sj}+\o_{sj}(X,Y)J_{is}].
\end{split}
\end{equation}
We take $i\ne j$, apply this to some vector $Z$ and take the scalar
product with $J_{ij}(W)$ to obtain
\begin{equation}\label{1r}
\begin{split}
R(X,Y,J_{ij}(Z),J_{ij}(W))-&R(X,Y,Z,W)=-2\o_{ij}(X,Y)g(J_{ij}(Z),W)\\+
&\sum_{s=1}^\rr[\o_{si}(X,Y)g(J_{si}(Z),W)+\o_{sj}(X,Y)g(J_{sj}(Z),W)].
\end{split}
\end{equation}

For $i\ne j$ we define the local two-forms $R^{ij}$ on $M$ by
\beq\label{not}R^{ij}(X,Y):=\sum_{a=1}^n
R(J_{ij}X_a,X_a,X,Y),\eeq
where $\{X_a\}$ denotes a local orthonormal frame on
$M$. In other words, $R^{ij}$ is twice the image of the 2-form $J_{ij}$ via
the curvature
endomorphism $R:\L^2M\to\L^2M$.
The first Bianchi identity easily shows that $R^{ij}(X,Y)=2\sum_{a=1}^n
R(X,X_a,J_{ij}X_a,Y)$.

(i) Assume that $\rr=4$. The image
$v:=\f(\o)$ of the volume element
$\o:=e_1\.e_2\.e_3\.e_4\in\Cl^0(E)$ is a parallel involution of $TM$
commuting with the $\Cl^0(E)$-action, so the tangent bundle of $M$
splits into a parallel direct sum $TM=T^+\oplus T^-$ of the $\pm 1$
eigen-distributions of $v$. By the de Rham decomposition theorem, $M$
is a Riemannian product $M=M^+\times M^-$. The restriction of $\f$ to
$\L^2_\pm E$ is trivial on $T^\pm$ and defines a rank 3 Clifford
structure on $M^\pm$. More explicitly, one can define a local
orthonormal frame
\beq\label{te}e ^\pm_1:=\frac12\bigg(e_1\wedge e_2\pm e_3\wedge
e_4\bigg),\quad
e^\pm_2:=\frac12\bigg(e_1\wedge e_3\mp e_2\wedge e_4\bigg),\quad
e ^\pm_3:=\frac12\bigg(e_1\wedge e_4\pm e_2\wedge e_3\bigg)\eeq
of  $\L^2_\pm E$ and it is clear that the local endomorphisms
$J^\pm_{ij}:=\f(e ^\pm_i)\circ\f(e ^\pm_j)$ vanish on $M^\pm$ and
satisfy the quaternionic relations on $M^\mp$. In fact it is
straightforward to check the relations
\beq\label{tj}J ^\pm_{12}=\pm\frac12\bigg(J_{14}\pm J_{23}\bigg),\quad
J^\pm_{31}=\pm\frac12\bigg(J_{13}\mp J_{24}\bigg),\quad
J^\pm_{23}=\pm\frac12\bigg(J_{12}\pm J_{34}\bigg)
\eeq
This shows that $M$ is a Riemannian product of two
quaternion-K\"ahler manifolds.

For later use, we remark that the curvature forms $\o ^\pm_{ij}$, $1\le
i,j\le 3$ of the
connection on $\L^2_\pm E$ with respect to the local frame $\{e_i^\pm\}$ are 
related to the forms $\o_{ij}$ by
\beq\label{o}\o ^\pm_{12}=\pm(\o_{14}\pm\o_{23}),\qquad
\o ^\pm_{31}=\pm(\o_{13}\mp\o_{24}),\qquad
\o ^\pm_{23}=\pm(\o_{12}\pm\o_{34}).
\eeq

(ii) Assume now that $\rr\ne 4$.
Let us choose some $k$ different from $i$ and $j$. Taking $Z=X_a$,
$W=J_{ik}(X_a)$, summing over $a$ in
\eqref{1r} and using \eqref{cstr} yields
\beq\label{2r}
2R^{ik}=\sum_{s=1}^\rr[\o_{si}<J_{si},J_{ik}>+\,\o_{sj}<J_{sj},J_{ik}>]\\
=n\o_{ik}.
\eeq
Taking now $Y=Z=X_a$ and summing over $a$ in
\eqref{1r} yields
$$\frac12R^{ij}(X,J_{ij}W)=\Ric(X,W)+2\o_{ij}(X,J_{ij}(W))-\sum_{s=1}^\rr[
\o_{si}(X,J_{si}(W))+\o_{sj}(X,J_{sj}(W))].$$
We identify 2-forms and endomorphisms on $M$ using $g$. The
previous relation reads
$$-\frac12 J_{ij}\circ R^{ij}=\Ric-2J_{ij}\circ \o_{ij}+\sum_{s=1}^\rr[
J_{si}\circ \o_{si}+J_{sj}\circ \o_{sj}],$$ so taking \eqref{2r}
into account we get for every $i\ne j$ \beq\label{3r}
0=\Ric+(n/4-2)J_{ij}\circ \o_{ij}+\sum_{s=1}^\rr[ J_{si}\circ
\o_{si}+J_{sj}\circ \o_{sj}]. \eeq It turns out that this system in
the unknown endomorphisms $J_{ij}\circ \o_{ij}$ has a unique solution for $n>8$.
Indeed, if we denote by $S_i:=\sum_{s=1}^\rr J_{si}\circ \o_{si}$
and sum over $j$ in \eqref{3r}, we get
$$0=\rr\Ric+ (n/4-2)S_i+\rr S_i+\sum_{j=1}^\rr S_j,$$
so $S_i=S_j$ for all $i,j$. From \eqref{3r} again we see that
$J_{ij}\circ \o_{ij}$ are all equal for $i\ne j$, and thus proportional
with $\Ric$:
\beq\label{4r}J_{ij}\circ \o_{ij}=\frac{1}{4-n/4-2\rr}\Ric,
\qquad\forall\ i\ne j.\eeq
Since the right term is symmetric, the two skew-symmetric endomorphisms
from the left term commute, so $J_{ij}$ commutes with $\Ric$ for all $i,j$.
This, in turn, implies like in Lemma \ref{or} above that
\beq\label{5r}<\o_{ij},J_{kl}>=0\
\mbox{unless } i=k\ne j=l\ \mbox{or}\ i=l\ne k=j.\eeq

We finally choose $k$ different from $i$ and $j$,
take $X=J_{ik}(X_a)$, $Y=X_a$, sum over $a$ in \eqref{1r} and use
\eqref{5r} to obtain
$$-J_{ij}\circ R^{ik}\circ J_{ij}-R^{ik}=-<\o_{ki},J_{ik}>J_{ki}.$$
By \eqref{2r} this reads
$$n\o_{ki}=-<\o_{ki},J_{ik}>J_{ki}$$
and \eqref{4r} then implies on the one hand that $M$ is Einstein
and on the other hand that the Ricci tensor does not vanish, since otherwise
$\n^E$ would be flat.

There exists thus a non-zero constant $\kk$ such that
\beq\label{h}\o_{ij}=\kk J_{ij}\eeq
for all $i\ne j$. This is equivalent to the statement (a).

We will now prove (iib) and (iii) simultaneously. From now on
$n$ might be equal to 8, but we assume that (a) holds.
We can re-express \eqref{0r}
and \eqref{4r} as
\begin{equation}\label{0r1}
R_{X,Y}J_{ij}=\kk\sum_{s=1}^\rr[g(J_{si}(X),Y)J_{sj}+g(J_{sj}(X),Y)J_{is}].
\end{equation}
and
\begin{equation}\label{0r2}
\Ric=\kk(n/4+2\rr-4).
\end{equation}
Assume that $M$ were reducible, i.e. that $TM$ is the direct sum of
two parallel distributions $T_1$ and $T_2$.
For all $X\in T_1$ and $Y\in T_2$ we have $R_{X,Y}=0$,
so \eqref{0r1} implies
$$0=\kk\sum_{s=1}^\rr[g(J_{si}(X),Y)J_{sj}+g(J_{sj}(X),Y)J_{is}].$$
Taking the scalar product with $J_{ik}$ for some $k\ne i,j$ and using
\eqref{5r} yields $$0=g(J_{kj}(X),Y).$$ This shows that each $J_{kj}$,
and hence the whole even Clifford structure, preserves
the splitting $TM=T_1\oplus T_2$.  In other words, each integral leaf $M_i$
of $T_i$ ($i=1,2$) carries a parallel
even Clifford structure. Notice that the relations $\o_{ij}=\kk
J_{ij}$ for all $i\ne j$ continue to hold on $M_1$ and $M_2$. Formula
\eqref{3r} then shows that the Ricci tensor of each factor
$T_i$ must satisfy $\Ric^{T_i}=\kk(\dim(T_i)/4+2\rr-4)$, which of course
contradicts \eqref{0r2}. This finishes the proof of (iib) and (iii).
\r

In order to proceed we need the following algebraic interpretation:

\begin{epr}\label{alg-int} Let $(M^n,g)$ be a simply connected
  Riemannian manifold with holonomy group $H:=\Hol(M)$ acting on
  $\RM^n$. A parallel rank $\rr$ ($3\le \rr\ne 4$) even Clifford
  structure on $M$ is 
  equivalent to an orthogonal representation $\rho:H\to \SO(\rr)$ of
  $H$ on $\RM^\rr$
  together with an $H$-equivariant algebra morphism
  $\ff:\Cl^0_\rr\to\End(\RM^n)$ mapping $\so(\rr)\subset \Cl^0_\rr$
  into $\so(n)\subset \End(\RM^n)$.
\end{epr}

\bp Assume that $\rho$ and $\ff$ satisfy the conditions above. Let $P$
be the holonomy bundle of $(M,g)$ through some orthonormal frame $u_0$,
with structure group $H$. The Levi-Civita connection of $M$ restricts to
$P$ and induces a connection on the Euclidean bundle
$E:=P\times_\rho \RM^\rr$. 
The bundle morphism
$$\f:\Cl^0(E)\to\End(TM),\qquad [u,a]\mapsto [u,\ff(a)]$$
is well-defined since $\ff$ is $H$-equivariant and clearly induces a
parallel rank $\rr$ even Clifford structure on $(M,g)$.

Conversely, if $(E,\n^E)$ defines a parallel even Clifford structure
on $M$, we claim that $E$ is associated to the holonomy bundle $P$
through $u_0$ and that $\n^E$ corresponds to the Levi-Civita
connection. Let $x_0$ be the base point of $u_0$, let $\G$ be the
based loop space at $x_0$ and let $\G_0$ be the kernel of the
holonomy morphism $\G\to H$. The parallel transport with respect to
$\n ^E$ of $E_{x_0}$ along curves in $\G$ defines a group morphism
$\tilde\rho:\G\to \SO(E_{x_0})$. If $\c\in\G_0$, the fact that
$(E,\n^E)$ is a parallel even Clifford structure is equivalent to
$\f(\L^2(\tilde\rho(\c))(\o))=\f(\o)$ for all $\o\in\L^2(E)$. Since
$\so(\rr)$ is simple for $3\le \rr\ne 4$, the map $\f$ is injective.
The relation above reduces to $\L^2(\tilde\rho(\c))=\id$, thus
to $\tilde\rho(\c)=\id$. This shows that $\G_0=\Ker(\tilde\rho)$, so
by taking the quotient, $\tilde\rho$ defines a faithful orthogonal
representation $\rho$ of $H=\G/\G_0$ on $E_{x_0}$. It is easy
to check that the map $P\times_\rho E_{x_0}\to E$ given
by
$$[u,e]\mapsto \tau ^E_\c(e),$$
where $\c$ is any curve in $M$ whose horizontal lift to $P$ through
$u_0$ ends at $u$ and $\tau ^E_\c$ denotes the parallel transport on
$E$ with respect to $\n^E$ along $\c$, is a well-defined bundle
morphism preserving the covariant derivatives. The existence of the
$H$-equivariant algebra morphism
  $\ff:\Cl^0_\rr\to\End(\RM^n)$ mapping $\so(\rr)\subset \Cl^0_\rr$
  into $\so(n)\subset \End(\RM^n)$ is now straightforward.
\r

It is easy to check that this result holds {\em verbatim} for
projective even Clifford structures, by replacing orthogonal
representations with projective ones. 
Notice that if $\rho:H\to \PSO(\rr)$ is a projective representation,
$\L^2\rho$ is a linear representation, so the vector bundle 
$\L^2E:=P\times_{\L^2\rho}\L^2(\RM^\rr)$ is globally defined, even
though $E:=P\times_\rho\RM^\rr$ is only locally defined.

\begin{ecor}\label{summand}
Assume that $(M^n,g)$ satisfies the hypotheses of Proposition $\ref{p37}$.
Then the Lie algebra $\hh$ of the holonomy group $H$ (associated to some
holonomy bundle $P$) is a direct sum of Lie sub-algebras, one of which
is isomorphic to $\so(\rr)$.
\end{ecor}
\begin{proof}
Every orthonormal frame $u_0\in P$ over $x_0\in M$, defines a natural
Lie algebra isomorphism from $\so(n)$ to $\L^2M_{x_0}$. In this way,
the holonomy algebra $\hh$ is naturally
identified with a sub-algebra of $\L^2M_{x_0}$ and the image
$\k$ of $\so(\rr)$ through the map $\ff$ defined in
Proposition \ref{alg-int} is naturally
identified with $\f(\L^2E_{x_0})$.

The
Ambrose-Singer Theorem (\cite{kn1}, Thm. 8.1 Ch.II) shows that $\hh$
contains the image of $\L^2M_{x_0}$ through the curvature
endomorphism. With the notation \eqref{not}, we thus get
$(R^{ij})_{x_0}\in\hh$ for all $1\le i,j\le \rr$. Taking
\eqref{2r} and \eqref{h} into account shows that $\k\subset \hh$.

Moreover, by Proposition \ref{alg-int},  $\k$ is an ideal of $\hh$.
Since $\hh$ is the Lie algebra of a compact Lie group, we immediately
obtain the Lie algebra decomposition $\hh=\k\oplus\k^\perp$, where
$\k^\perp$ is the orthogonal complement of $\k$ in $\hh$ with respect to any
$\ad_{\hh}$-invariant metric on $\hh$.
\r

We are now ready for the first important result of this section.

\begin{ath}
A Riemannian manifold $(M^n,g)$ carrying a parallel
non-flat even Clifford structure $(E,\n^E)$ of rank $\rr\ge 5$ is
either locally symmetric or $8$-dimensional.
\end{ath}
\bp
Assume that $M$ is not locally symmetric. By replacing $M$ with its
universal cover, we may assume that $M$ is simply connected. 
According to Proposition \ref{p37},
$M$ has irreducible holonomy and non-vanishing scalar curvature. The
Berger-Simons Holonomy Theorem implies that there are exactly three
possibilities for the holonomy group $H$ of $M$: $H=\SO(n),\
H=\U(n/2)$ or $H=\Sp(n/4)\.\Sp(1)$. The second exterior power of the
holonomy representation is of course irreducible in the first case and
decomposes as
$$\so(n)=\su(n/2)\oplus \RM\oplus \pp_1,$$
$$\so(n)=\sp(n/4)\oplus\sp(1)\oplus \pp_2$$
in the latter two cases. A summand isomorphic to some $\so(\rr)$
($\rr\ge 5$) occurs in the above decompositions if and only if $\rr=n$
in the first case, or is obtained from
the low-dimensional isomorphisms
$$\su(n/2)\simeq \so(\rr)\qquad\hbox{for }n=8 \hbox{ and }
\rr=6,$$
$$\sp(n/4)\simeq\so(\rr)\qquad\hbox{for }n=8\hbox{ and }\rr=5.$$
In the latter cases one has $n=8$, so we are left with the case when
$M$ has generic holonomy $\SO(n)$. By
Proposition \ref{alg-int}, $\RM^n$ inherits a $\Cl^0_n$-module
structure, which for dimensional reasons may only occur when $n=8$.
\r

Using this result we will now obtain the classification of
complete simply connected 
manifolds with parallel rank $\rr$ even Clifford structures. From the
above discussion it is enough to consider the cases when $\rr\ge 5$
and either $\dim(M)=8$ or $M$ is symmetric.

{\bf Case 1. $\dim(M)=8$.} Proposition \ref{alg-int} has several
consequences:
\begin{itemize}
\item[(a)] $\RM^8$ is a
$\Cl^0_\rr$ representation, thus  $\;5\le\rr\le 8$.
\item[(b)] The inclusion
$\ff:\so(\rr)\to\so(8)$ is defined by the spin (or half-spin for
$\rr=8$) representation.
\item[(c)] The holonomy group $H$ is contained in
the connected component of the identity, called $N^0_{\SO(8)}\so(\rr)$, of the
normalizer of $\so(\rr)$ in $\SO(8)$, acting on its Lie algebra by
the adjoint representation.
\end{itemize}
Using again the low-dimensional
isomorphisms $\so(5)\simeq \sp(2)$ and $\so(6)\simeq\su(4)$ we easily
get
$$N^0_{\SO(8)}\so(5)=\Sp(2)\.\Sp(1),\qquad N^0_{\SO(8)}\so(6)=\U(4),
\qquad N^0_{\SO(8)}\so(7)=\Spin(7).$$
Thus a necessary condition for a simply connected 
8-dimensional manifold to carry a
parallel even Clifford structure of rank $\rr$ is that $M$ is
quaternion-K\"ahler for $\rr=5$, K\"ahler for $\rr=6$ and has holonomy
contained in $\Spin(7)$ for $\rr=7$ (no condition at all for $\rr=8$).
Conversely, if $M$ satisfies one of these conditions for $\rr
=5,6,7$ or is an arbitrary manifold in the case $\rr = 8$, we define
$E$ to be associated to the holonomy bundle of $M$ with respect to
the following representations of the holonomy group:
$$
\begin{array}{rl}
\rr = 5: & \Sp(2)\.\Sp(1)\to \SO(5),\ a\.b\mapsto \xi(a),
\;\mbox{where}\quad
  \xi:\Sp(2)\simeq\Spin(5)\to\SO(5)\\
  &  \mbox{is the spin covering.} \\
\rr = 6: &  \U(4)\to\PSO(6) \;\mbox{induced by taking the} \; \ZM_4
\;\mbox{quotient in the
  projection}\\
  & \mbox{onto the first factor in} \quad
  \SU(4)\times\U(1)\to\SU(4)\simeq \Spin(5).\\ 
\rr = 7: & \mbox{The spin covering} \quad \Spin(7)\to\SO(7).\\
\rr = 8: & \mbox{One of the two representations}\quad
\SO(8)\to\PSO(8) \;\mbox{obtained by taking the} \\
&\;\ZM_2\; \mbox{quotient in the half-spin representations}\;
  \Spin(8)\to\SO(\Delta_\pm).
\end{array}
$$

Notice that for $\rr=6$ and $\rr=8$ the defining representation of
$E$ is projective, so $E$ is only locally defined if $M$ is
non-spin. On the contrary, if $M$ is spin then $E$ is a well-defined
vector bundle, associated to the spin holonomy bundle of $M$.

The attentive reader might have noticed the subtlety of the case
$\rr=8$. In all other cases the equivariant Lie algebra morphism $\ff$
is constructed by identifying $\so(\rr)$ with a factor of the Lie
algebra of the holonomy group acting on $\RM^8$ by the {\em spin}
representation (therefore extending to a representation of the even
Clifford algebra). For $\rr=8$ however, the holonomy representation
{\em is not} the spin representation. What still makes things work in
this case is the {\em triality} of the $\so(8)$ representations, which
is an outer automorphism of $\Spin(8)$ interchanging its three
non-equivalent representations on $\RM^8$. In this way, on a
8-dimensional spin manifold one has six Clifford actions: The
Clifford algebra bundle of $TM$ acts on the half spinor bundles
$\Sigma_\pm M$,  
$\Cl^0(\Sigma_+M)$ acts on $TM$ and $\Sigma_-M$, and 
$\Cl^0(\Sigma_-M)$ acts on $TM$ and $\Sigma_+M$. Of course, when $M$
is not spin, 
among the six Clifford actions above, only the third and the fifth ones
are globally defined.

According to Proposition \ref{alg-int}, the argument above can
be expressed as
follows: We denote by $\xi:\Spin(8)\to\SO(8)$ the spin covering and
by $\delta^\pm:\Spin(8)\to\SO(8)$ the half-spin representations. If
$H\subset\SO(8)$ is the holonomy group of $M$, let
$\rho:H\to\PSO(8)$ denote the restriction to $H$ of the $\ZM_2$-quotient
of $\d^+$. The isomorphism $\ff:\so(8)\to\so(8)$,
$\ff=\xi_*\circ(\delta^+_*)^{-1}$ is tautologically equivariant with
respect to the representations of $H$ on $\so(8)$ induced by
$\rho$ and $\xi$ respectively, and it extends to a Clifford action
due to triality.

{\bf Case 2. $M=G/H$ is symmetric.}
According to Proposition \ref{alg-int} and Corollary \ref{summand}, 
there are two necessary
conditions for $M$ to carry a parallel even Clifford structure of
rank $\rr\ge 5$:
\begin{itemize}
\item[(a)]  $\so(\rr)$ occurs as a summand in the Lie
algebra $\hh$ of the isotropy group $H$.
\item[(b)] The dimension of $M$ has to be a multiple of
the dimension $N_0(\rr)$ of the irreducible $\Cl^0_\rr$
representation.
\end{itemize}
Notice that Proposition \ref{alg-int} shows that if
$M=G/H$ is a compact symmetric space solution of our problem, its
non-compact dual $G^*/H$ is a solution too, since the isotropy
representations are the same. We will thus investigate only the
symmetric spaces of compact type.

After a cross-check in the tables of symmetric spaces of Type I and II
(\cite{besse}, pp. 312-317) we are left with the following cases:

(1) $G=\SU(n)$, $H=\SO(n)$. Condition (a) is verified for $\rr=n$
  but it is easy to check that $\dim(M)=(\rr-1)(\rr+2)/2$ cannot be a
  multiple of $N_0(\rr)$.
  
(2) $G=\SU(2n)$, $H=\Sp(n)$. Condition (a) is verified for $n=2$ and
  $\rr=5$, but $\dim(M)=5$ is not a multiple of $N_0(5)=8$.
 
(3) $G=\SU(p+q)$, $H=S(\U(p)\times \U(q))$. Both conditions are
  verified for $p=4$, $\rr=6$ and arbitrary $q$.
 
(4) $G=\SO(p+q)$, $H=\SO(p)\times \SO(q)$. By condition (a) one can
  assume $\rr=p\ge 5$. The isotropy representation is the tensor
  product $\RM^{pq}$ of the standard representations of $\SO(p)$ and
  $\SO(q)$.
Assume that $p\ne 8$. It is well known that the group $\SO(p)$ has
exactly one non-trivial representation on $\RM^p$. This is due to the
fact that $\SO(p)$ has no outer automorphisms for $p$ odd, while for
$p$ even the only outer automorphisms are the conjugations by matrices
in $\mathrm{O}(p)\setminus \SO(p)$. Restricting our attention to the
subgroup
$\SO(p)$ of the holonomy group $H$, the map $\ff$ given by Proposition
\ref{alg-int}
defines an $\SO(p)$-equivariant representation of $\so(p)$ on
$\RM^p\oplus\ldots\oplus\RM^p$ ($q$ times) and is thus defined by $q^2$
equivariant components $\ff_{ij}:\so(p)\to \End(\RM^p)$. It is easy to
see that each $\ff_{ij}$ is then scalar: $\ff_{ij}(A)=\l_{ij}A$ for
all $A\in\so(p)$. Finally, the fact that $\ff$ extends to the Clifford
algebra implies that $\ff(A)^2=-\id$ for $A=\xi_*(e_1\.e_2)$ (here
$\xi$ denotes the spin covering $\Spin(p)\to\SO(p)$), and this is
impossible since
$$(\ff(A)^2)_{ij}=\sum_{k=1}^q\l_{ik}\l_{kj}A^2,$$
and $A^2$ is not a multiple of the identity. Thus $\rr=p=8$
is the only admissible case.
 
(5) $G=\SO(2n)$, $H=\U(n)$. Condition (a) is verified for $n=4$ and
  $\rr=6$, but $\dim(M)=12$ is not a multiple of $N_0(6)=8$.
 
(6) $G=\Sp(n)$, $H=\U(n)$. Condition (a) is verified for $n=4$ and
  $\rr=6$, but $\dim(M)=20$ is not a multiple of $N_0(6)=8$.
 
(7) $G=\Sp(p+q)$, $H=\Sp(p)\times \Sp(q)$. Both conditions are
  verified for $p=2$, $\rr=5$ and arbitrary $q$.
 
(8) If $G$ is one of the exceptional simple Lie groups $\F_4$, $\E_6$,
  $\E_7$, $\E_8$, both conditions are simultaneously verified for
  $H=\Spin(9)$, $\Spin(10)\times \U(1)$, $\Spin(12)\times \SU(2)$ and
  $\Spin(16)$ respectively. The corresponding symmetric spaces are
  exactly Rosenfeld's elliptic projective planes $\Ca \PM^2$,
  $(\CM\otimes
\Ca) \PM^2$, $(\HM\otimes \Ca) \PM^2$ and $(\Ca\otimes \Ca) \PM^2$.
 
(9) Finally, no symmetric space of type II (i.e. $M=H\times H/H$)
  can occur: condition (a) is satisfied for $H=\SU(4)$, $\rr=6$ and
  $H=\SO(n)$, $\rr=n$ but the dimension of $M$ is $15$ in the first
  case and $n(n-1)/2$ in the second case, so condition (b) does not
  hold.

The only candidates of symmetric spaces carrying parallel even Clifford
structures of rank $\rr\ge 5$ are thus those from cases (3), (4), (7)
and (8). Conversely, all these spaces carry a (projective) parallel
even Clifford 
structure. This is due to the fact that the restriction of the
infinitesimal isotropy representation to the $\so(\rr)$ summand is the
spin representation in all cases except for $\so(8)$, where the triality
argument applies. Summarizing, we have
proved the following
\begin{ath}\label{class}
The list of complete simply connected Riemannian manifolds $M$ carrying a
parallel rank $\rr$ even Clifford structure is given in the tables
below.
\end{ath}

\begin{center}
\begin{tabular}{|r|l|c|}\hline
$\rr$    &  $M$   &  dimension of $M$ \\
 \hline\hline
2            &   K\"ahler  &  $2m,\ m\ge 1$      \\
\hline
3 and 4           &  hyper-K\"ahler  &  $4q,\ q\ge 1$      \\
\hline
4            &  reducible hyper-K\"ahler &
$4(q^++q^-),$   $ q^+\ge 1$, $q^-\ge 1$   \\
\hline
arbitrary          & $\Cl^0_{\rr}$ representation space    &
multiple of $N_0(\rr)$     \\

\hline
\end{tabular}
\vskip .2cm
 Table 1. Manifolds with a flat even Clifford
  structure
\end{center}
\vs

\begin{center}
\begin{tabular}{|r|l|l|c|}\hline
$\rr$    & type of $E$  &  $M$   &  dimension of $M$ \\
 \hline\hline
2         &      & K\"ahler  &  $2m,\ m\ge 1$      \\
\hline
3         & projective if $M\ne \HM\PM^q$     &  quaternion-K\"ahler 
(QK)  &  $4q,\ q\ge 1$      \\
\hline
4          & projective if $M\ne \HM\PM^{q^+}\times\HM\PM^{q^-}$  &  product
of two QK manifolds& 
$4(q^++q^-)$      \\
\hline\hline
5          &     &  QK    &  8      \\
\hline
6          &projective if $M$ non-spin   &   K\"ahler  &  8      \\
\hline
7          &     &  $\Spin(7)$ holonomy   &  8      \\
\hline
8          & projective if $M$ non-spin   &  Riemannian  &  8      \\
\hline\hline
5   &  & $\Sp(k+2)/\Sp(k)\times\Sp(2)$ & $8k,\ k\ge 2$ \\
\hline
6   & projective  & $\SU(k+4)/{\rm S}(\U(k)\times\U(4))$ & $8k,\ k\ge 2$ \\
\hline
8   & projective if $k$ odd & $\SO(k+8)/\SO(k)\times\SO(8)$ & $8k,\ k\ge 2$ \\
\hline\hline
9    &   & $\hskip1.15cm\Ca \PM^2=\F_4/\Spin(9) $  &  16 \\
\hline
10   &   & $(\CM\otimes \Ca) \PM^2=\E_6/\Spin(10)\.\U(1)$   &  32 \\
\hline
12   &   & $(\HM\otimes \Ca) \PM^2=\E_7/\Spin(12)\.\SU(2)$   &  64 \\
\hline
16   &   & $(\Ca\otimes \Ca) \PM^2=\E_8/\Spin^+(16)$   &  128 \\
\hline
\end{tabular}

Table 2. Manifolds with a parallel non-flat even Clifford
  structure \footnote{In this table we adopt the convention that the
 QK condition is empty in dimension $4$. For the sake of simplicity,
 we have omitted in Table~2 the 
    non-compact duals of the compact symmetric spaces. The
    meticulous reader should add the spaces obtained by replacing
    $\Sp(k+8)$, $\SU(k+4)$, $\SO(k+8)$, $\F_4$, $\E_6$,  $\E_7$ and $\E_8$
    in the last seven rows
    with  $\Sp(k,8)$, $\SU(k,4)$, $\SO_0(k,8)$, $\F_4^{-20}$,
    $\E_6^{-14}$,  $\E_7^{-5}$ and $\E_8^8$ respectively.}
\end{center}
\vs 

We end up this section with the classification of manifolds carrying
parallel Clifford structures.

\begin{ath}\label{t23}
A simply connected Riemannian manifold $(M^n,g)$ carries a parallel
rank $\rr$ Clifford structure if and only if one of the following
(non-exclusive) cases occurs:
\begin{enumerate}
\item $\rr=1$ and $M$ is K\"ahler.
\item $\rr=2$ and either $n=4$ and $M$ is K\"ahler or $n\ge 8$ and $M$
  is hyper-K\"ahler.
\item $\rr=3$ and $M$ is quaternion-K\"ahler.
\item $\rr=4$, $n=8$ and $M$ is a product of two Ricci-flat K\"ahler surfaces.
\item $\rr=5$, $n=8$ and $M$ is hyper-K\"ahler.
\item $\rr=6$, $n=8$ and $M$ is K\"ahler Ricci-flat.
\item $\rr=7$ and $M$ is an $8$-dimensional manifold with $\Spin(7)$ holonomy.
\item $\rr$ is arbitrary and $M$ is flat, isometric to a
  representation of the Clifford algebra $\Cl_\rr$.
\end{enumerate}
\end{ath}

\begin{proof} Assume that $(M^n,g)$
  carries a rank $\rr$ parallel Clifford structure $(E,h)\subset
  (\L^2M,\frac1n g)$. The image by $\f:\Cl(E,h)\to\End(TM)$ of the
  volume element is 
  a parallel endomorphism $v$ of $TM$ which satisfies $v\circ
  v=(-1)^{\frac{\rr(\rr+1)}{2}}$ and commutes (resp. anti-commutes)
  with every element of $E$ for $\rr$ odd (resp. even). We start by
  considering the cases $\rr\le 4$.

\noindent$\bullet$ {$\rr=1$.} It was already noticed that a parallel
rank 1 Clifford structure
corresponds to a K\"ahler structure on $M$.

\noindent$\bullet$ {$\rr=2$.} The rank 2 Clifford structure $E$
induces a rank 3 Clifford structure $E ':=E\oplus\L^2E$ on $M$.
Explicitly, if $\{e_1,e_2\}$ is a local orthonormal basis of $E$,
then $e_3:=e_1\circ e_2$ is independent of the chosen basis and
$\{e_1,e_2,e_3\}$ satisfy the quaternionic relations. Moreover,
$e_3=v$ is a parallel endomorphism of $TM$, so $(M,g)$ is K\"ahler.
In the notation of Proposition \ref{p37} we have
$\o_{13}=\o_{23}=0$. Formula \eqref{3r} yields
$$
0=\Ric+n/4J_{ij}\circ \o_{ij}+
J_{si}\circ \o_{si}+J_{sj}\circ \o_{sj}
$$
for every permutation $\{i,j,s\}$ of $\{1,2,3\}$. If $n>4$ this system
shows that $\o_{12}=0$, so $M$ is hyper-K\"ahler. Conversely, if either
$n=4$ and $(M,g,J)$ is K\"ahler, or $n>4$ and $(M,g,I,J,K)$ is
hyper-K\"ahler, then $E=\L^{(2,0)+(0,2)}M$ in the first case, or
$E=<I,K>$ in the second case, define a rank 2 parallel Clifford
structure on $M$.

\noindent$\bullet$ {$\rr=3$.} It was already noticed that because
of the isomorphism $\Lambda^2E \cong E$,
every rank 3 even Clifford structure is automatically a Clifford
structure, and corresponds to a quaternion-K\"ahler structure
(which, we recall, is an empty condition for $n=4$).

\noindent$\bullet$ $\rr=4$. The endomorphism $v$ is now a parallel
involution of $TM$ anti-commuting with every element of the Clifford bundle
$E\subset \L^2M$. Correspondingly, the tangent bundle of $M$ splits in
a parallel
direct sum $TM=T^+\oplus T^-$, such that $v\res_{T^\pm}=\pm\id$. If
  we denote by $J_i$, $1\le i\le 4$ a
  local orthonormal basis of $E$, each $J_i$ maps $T^\pm$ to
  $T^\mp$. The de Rham decomposition theorem shows that $M$ is a
  Riemannian product $M=M^+\times M^-$ and $TM^\pm= T^\pm$. The
  Riemannian curvature tensor of $M$ is the sum of the two curvature
  tensors of $M^+$ and $M^-$: $R=R^++R^-$.
Let $\o_{ij}$ denote the curvature forms (with respect to the local
frame $\{J_i\}$) of the Levi-Civita
connection on $E$:
$$R_{X,Y}J_i=\sum_{j=1}^4\o_{ji}(X,Y)J_j.$$
We take $X,Y\in T^+$
and apply the previous relation to some $Z\in T^+$ and obtain
$$R^+_{X,Y}Z=\sum_{j=1}^4\o_{ji}(X,Y)J_jJ_iZ,\qquad \forall\ 1\le i\le 4\.$$
For $1\le j\le 3$ we denote by $\o_j:=\o_{j4}$ and
$I_j:=-J_jJ_4$.
Since by definition $v=J_1J_2J_3J_4$, it is easy to
check that  $I_j$ are anti-commuting almost complex
structures on $M^+$ satisfying the quaternionic relations $I_1I_2=I_3$
etc. The previous curvature relation reads
\beq\label{k}R^+_{X,Y}=-\sum_{j=1}^3\o_{j}(X,Y)I_j,\qquad \forall\
X,Y\in T^+.
\eeq
The symmetry by pairs of $R^+$ implies that
$\o_i=\sum_{j=1}^3a_{ji}I_j$ for some smooth functions $a_{ij}$
satisfying $a_{ij}=a_{ji}$. Moreover, the first Bianchi identity
applied to \eqref{k} yields
\beq\label{su}\sum_{i,j=1}^3a_{ij} I_i\wedge I_j=0.
\eeq
If $\dim(M^+)>4$, we may choose non-vanishing vectors $X,Y\in T^+$
such that $Y$ is
orthogonal to $X$ and to $I_iX$ for $i=1,2,3$. Applying \eqref{su} to
$X,I_i X, Y, I_jY$ yields $a_{ij}+a_{ji}=0$, so $\o_j=0$. By \eqref{k}
we get $R^+=0$ and similarly $R^-=0$, so $M$ is flat. It remains to
study the case $\dim(M^+)=4$. In this case $I_1,I_2$ and $I_3$ are a
basis of the space of self-dual 2-forms $\L^2_+M^+$, so \eqref{k} is
equivalent to the fact that $M^+$ is self-dual and has vanishing
Ricci tensor (see e.g. \cite{besse}, p.51). In other words, $M^+$ is
K\"ahler (with respect to any parallel 2-form in $\L^2_-M^+$) and
Ricci-flat, and the same holds of course for $M^-$.

Conversely, assume that $M=M^+\times M^-$ is a Riemannian product of
two simply connected Ricci-flat K\"ahler surfaces. The holonomy of $M$
is then a subgroup of $\SU(2)\times \SU(2)\simeq\Spin(4)$, so the
frame bundle of $M$ and the
Levi-Civita connection reduce to a principal $\SU(2)\times
\SU(2)$-bundle $P$. Let $\xi$
denote the representation of $\Spin(4)$ on
  $\RM^4$ coming from the spin covering
  $\Spin(4)\to\SO(4)$ and let $\rho$ denote the representation of
  $\Spin(4)$ on $\so(8)$ obtained by
  restricting the adjoint action of $\SO(8)$ to
  $\Spin(4)\simeq\SU(2)\times \SU(2)\subset \SO(8)$.
The irreducible representation of $\Cl_4$ on $\RM^8$
  defines a $\Spin(4)$-equivariant map from $\RM^4$ to $\so(8)$ (with
  respect to the above actions of $\Spin(4)$). The above map defines
  an embedding
  of the rank 4 vector bundle $E:=P\times _{\xi}\RM^4$ into $\L^2 M=
P\times _{\rho}\so(8)$, which is by construction a parallel
Clifford structure on $M$.

For $\rr\ge 5$ we will use the fact that $E$ defines tautologically a
  rank $\rr$ parallel even Clifford structure on $M$, and apply
  Theorems \ref{t35} and \ref{class} to reduce the study to manifolds
  appearing in Table 2.

\noindent$\bullet$ $\rr=5$. The volume element $v$ defines a K\"ahler structure
  on $M$ in this case. The quaternionic Grassmannians
  $\Sp(k+8)/\Sp(k)\.\Sp(2)$ are
  obviously not K\"ahler (since the Lie algebra of the isometry group 
  of every K\"ahler symmetric space has a non-trivial center), so it 
  remains to examine the case $n=8$, when,
  according to Theorem \ref{class}, $M$ is quaternion-K\"ahler. More
  explicitly, if $E$ is the rank 5 Clifford
  bundle, $\f(\L^2 E)$ is a Lie sub-algebra of $\End^-(TM)\simeq \so(8)$
  isomorphic to $\so(5)\simeq\sp(2)$ and its centralizer is a Lie
  sub-algebra $\mathfrak s$ of $\End^-(TM)$ isomorphic to $\so(3)$, defining
  a quaternion-K\"ahler structure. Moreover $v$ belongs to
  $\mathfrak s$ (being the image of a central element in the Clifford
  algebra bundle of $E$), so we easily see that its orthogonal complement
  $v^\perp$ in $\mathfrak s$ defines a rank 2 parallel Clifford
  structure on $M$. By the case $\rr=2$ above, $M$ is then
  hyper-K\"ahler.

Conversely, every 8-dimensional hyper-K\"ahler
  manifold carries parallel Clifford structures of rank 5 obtained as
  follows. Let $\xi$ denote the representation of $\Spin(5)$ on
  $\RM^5$ coming from the spin covering
  $\Spin(5)\to\SO(5)$ and let $\rho$ denote the representation of
  $\Spin(5)$ on $\so(8)$ obtained by
  restricting the adjoint action of $\SO(8)$ to
  $\Spin(5)\simeq\Sp(2)\subset \SO(8)$.
The irreducible representation of $\Cl_5$ on $\RM^8$
  defines a $\Spin(5)$-equivariant map from $\RM^5$ to $\so(8)$ (with
  respect to the above actions of $\Spin(5)$). If $P$ denotes the
  holonomy bundle of $M$ with structure group $\Sp(2)\simeq
  \Spin(5)$, the above map defines an embedding
  of the rank 5 vector bundle $E:=P\times _{\xi}\RM^5$ into $\L^2 M=
P\times _{\rho}\so(8)$, which is by construction a parallel
Clifford structure on $M$.

\noindent$\bullet$ $\rr=6$. The volume element $v$ is now a K\"ahler structure
  anti-commuting with every element of the Clifford bundle $E$. If
  we denote by $J_i$, $1\le i\le 6$ a
  local orthonormal basis of $E$, each $J_i$ is a
  2-form of type $(2,0)+(0,2)$ with respect to $v$, so the curvature
  endomorphism vanishes on $J_i$:
\beq\label{v}
0=R(J_i)(X,Y)=\sum_{a=1}^n R(J_iX_a,X_a,X,Y)=2\sum_{a=1}^n R(X,X_a,J_iX_a,Y).
\eeq
Let $\o_{ij}$ denote the curvature forms (with respect to the local
frame $\{J_i\}$) of the Levi-Civita
connection on $E$:
$$R_{X,Y}J_i=\sum_{j=1}^6\o_{ji}(X,Y)J_j.$$
We can express this as follows:
$$R(X,Y,J_iZ,J_iW)-R(X,Y,Z,W)=\sum_{j=1}^6\o_{ji}(X,Y)g(J_jZ,J_iW).$$
Taking the trace in $Y$ and $Z$ and using \eqref{v} yields
$$\Ric=-\sum_{j=1}^6J_j\circ J_i\circ\o_{ji}.$$
This relation, together with \eqref{3r}, shows that $\Ric=0$.

Conversely, every 8-dimensional Ricci-flat K\"ahler manifold carries
parallel Clifford structures of
  rank 6 defined by the $\Spin(6)\simeq\SU(4)$-equivariant embedding of $\RM^6$
  into $\so(8)$ coming from the irreducible representation of
  $\Cl_6$ on $\RM^8$, like in the case $\rr=5$.

\noindent$\bullet$ $\rr=7$. Theorem \ref{class} shows that  $M$ has to be an
  8-dimensional manifold with holonomy $\Spin(7)$. By an argument
  similar to the previous ones,
every such manifold carries parallel Clifford structures of
  rank 7 defined by the $\Spin(7)$-equivariant embedding of $\RM^7$
  into $\so(8)$ coming from one of the irreducible representations of
  $\Cl_7$ on $\RM^8$.

\noindent$\bullet$ $\rr=8$. The dimension of $M$ has to be at least
equal to 16 in this case (since the
  dimension of the irreducible $\Cl_8$-representation is 16). Moreover, the
  volume element $v$ is a parallel involution of
  $TM$ anti-commuting with every element of $E$, so $TM$ splits in a
  parallel direct sum of the $\pm 1$ eigen-distributions of $v$. This
  contradicts Proposition \ref{p37}.

\noindent$\bullet$ Finally, for $\rr\ge 9$, the spaces
appearing in the last four rows of Table 2 cannot carry a Clifford
structure since the dimension of the irreducible representation of
  $\Cl_\rr$ for $\rr=9,\ 10,\ 12,\ 16$ is $32,\ 64,\ 128,\ 256$
  respectively, which is exactly twice the dimension of the
  corresponding tangent spaces in each case.
\r

\vs

\section{Bundle-like curvature constancy}\label{scc}

As an application of Theorem \ref{class}, 
we classify in this section bundle-like metrics with curvature
constancy. We first 
show in Subsection \ref{rs} that every Riemannian submersion 
$Z\to M$ with totally geodesic fibres is associated to a $G$-principal
bundle $P\to 
M$ (where $G$ is the isometry group of some given fibre), which
carries a canonical $G$-invariant connection. The curvature
of this connection is a 2-form $\o$ on $M$ with values in the
adjoint bundle $\ad(P)$. We then compute the different components of
the Riemannian curvature tensor of $Z$ in terms of the Riemannian
curvature of $M$ and of the curvature form $\o$. 

Most of this material can be found in the literature (cf. \cite{GW},
see also \cite{ziller}), but we include it here for summing up the
notations, conventions and usual normalizations. Readers familiar with
Riemannian geometry can pass directly to Subsection \ref{cco}, where
we interpret the curvature constancy condition \eqref{cn} by the fact
that $\o$ defines a {\em parallel even Clifford structure} on $M$. The
classification is obtained in Subsection \ref{cla} by a case-by-case
analysis through the manifolds in Table 2.

\subsection{Riemannian submersions with totally geodesic
  fibres}\label{rs} 

Let $\pi:Z^{k+n}\to M^n$ be a Riemannian submersion with totally geodesic
fibres. Assume that $Z$ is complete. We denote by $Z_x:=\pi ^{-1}(x)$
the fibre of
$\pi$ over $x\in M$. From Theorem~1 in \cite{her}, all
fibres are isometric to some fixed Riemannian manifold $(F,g_F)$ and
$\pi$ is a locally trivial fibration with structure group the Lie
group $G:=\Iso(F)$ of isometries of $F$.

For every tangent vector $X\in T_xM$ and $z\in Z_x$, we denote
by $X^*$ its horizontal lift at $z$. For every curve $\c$
on $M$ and $z\in Z_{\c(0)}$ there exists a unique curve $\tilde
\c$ with $\tilde\c(0)=z$ whose tangent vector at $t$ is the
horizontal lift of $\dot{\c(t)}$ at $\tilde\c(t)$ for every $t$. This
is called the {\em horizontal lift} of $\c$ through $z$.
Hermann's result in \cite{her} mainly says that for every curve $\c$
on $M$, the mapping $\tau_t:Z_{\c(0)}\to Z_{\c(t)}$,
which maps $z$ to the value at $t$ of the horizontal lift of $\c$
through $z$, is an isometry between the two fibres, (each endowed with
the induced Riemannian metric).

We define the $G$-principal fibre bundle $P$ over $M$ as the set of isometries
from $F$ to the fibres of $\pi$:
$$
P:=\{u:F\to Z\ |\ \exists\  x\in M\ \mbox{such that}\ u\ \mbox{maps
  $F$ isometrically onto $Z_x$}\}.
$$
We denote by $p:P\to M$ the natural projection and by $P_x$ the fiber
of $p$ over $x$:
$$P_x:=\{u:F\to Z_x\ |\ u\ \mbox{is an isometry}\}.$$
The right action of $G=\Iso(F)$ on $P$ is given by $ua:=u\circ a$ for
every $u\in P$ and $a\in G$.

\begin{epr}\label{conn} {\em (Cf. \cite{GW}, Theorem 2.7.2)} The
  horizontal distribution on $Z$ induces a $G$-invariant connection on $P$.
\end{epr}
\begin{proof}

For $X\in T_xM$ and $u\in P_x$, we define its {\em horizontal lift}
$\t X\in T_uP$ as follows. Take any curve $x_t$ in $M$ such that
$X=\dot x_0$. The isometry $\tau_t$ between $Z_{x_0}$ and
$Z_{x_t}$ described above, defines a curve $u_t:=\tau_t\circ u$ which
obviously satisfies $p(u_t)=x_t$. We then set $\t X:=\dot u_0$ and
claim that this does not depend on the curve $x_t$. This is actually a
direct consequence of the following more general result:

\begin{elem}\label{ind}
Let $p:P\to M$ be a $G$-principal fibre bundle and assume that $G$ acts
effectively on some manifold $F$. Define $Z:=P\times_G F$ and for each
$f\in F$, the
smooth map $R_f:P\to Z$, $R_f(u)=u(f)$. Then a tangent vector $X\in
T_uP$ vanishes if and only if $p_*(X)=0$ and $(R_f)_*(X)=0$ for every
$f\in F$.
\end{elem}
\begin{proof} Since the result is local, one
  may assume that $P=M\times G$ is trivial and $u=(x,1)$. One can
  write $X=(X',X'')$,
  with $X'\in T_xM$ and $X''\in \gg$. Since $p_*(X)=0$, we get
  $X'=0$. From $(R_f)_*(X)=0$ we obtain $\exp(tX'')(f)=f$ for
  every $t\in \RM$ and $f\in F$.  If
  $X''$ were not zero, this would contradict the effectiveness of the
  action of $G$.
\r

Returning to our argument, we see that $p_*(\t X)=X$ and
$$(R_f)_*(\t X)=\frac{\del}{\del t}\Big|_{t=0}(u_t(f))=\frac{\del}{\del
  t}\Big|_{t=0}\tau_t(u(f))=X^*_{u(f)}$$
only depend on $X$, not on $x_t$. The map $T_xM\to T_uP$, $X\mapsto \t
X$ is thus well-defined for every $x\in M$ and $u\in p^{-1}(x)$. We
denote by $H_u$ the image of this map.

Lemma \ref{ind} also shows that $H_u$ is a vector subspace of $T_uP$,
supplementary to the tangent space to the fibre of $P$ through
$u$. The collection $\{H_u,u\in P\}$  is called the horizontal
distribution, and it is easy to see that it
is invariant under the action of $G$: If $a\in G$ $u\in P$ and $x_t$ is a
curve in $M$ with $x_0=p(u)$, then (denoting $X:=\dot x_0$):
$$(R_a)_*(\t X_u)=(R_a)_*\frac{\del}{\del t}\Big|_{t=0}(\tau_t \circ u)=
\frac{\del}{\del t}\Big|_{t=0}(\tau_t \circ ua)=\t X_{ua}.$$
This proves the proposition.
\r

We will now express the Riemannian curvature of $Z$ in terms of the
curvature of the connection on $P$ defined above (we will denote this
connection by $\theta$ in the sequel). In order to do this, we need to
introduce some notation. The adjoint bundle $\ad(P)$ of $P$, is the
vector bundle associated to $P$ via the adjoint representation of $G$
on its Lie algebra:
$$
 \ad(P):=P\times_{\ad}\gg,
$$
where for every $g\in G$, $\ad_g:\gg\to\gg$ is the differential at the
identity of $\Ad_g:G\to G$ defined as usually by
$\Ad_g(h):=ghg^{-1}$.
The curvature of the connection $\theta$ defined by Proposition
\ref{conn} is a $G$-equivariant 2-form $\t\o$ on $P$ with values in
$\gg$ or, equivalently, a 2-form $\o$ on $M$ with values in the vector bundle
$\ad(P)$, i.e. a section  of $\Lambda^2M\otimes
{\ad(P)}$.
The forms $\o$ and $\t\o$
are related by
\beq\label{om}[u, \t\o(\t X,\t Y)_u]=\o(X, Y)_{p(u)},
\eeq
where $X,Y \in T_{p(u)}M$ are tangent vectors on $M$ with horizontal lifts
$\t X, \t Y \in T_uP$ to tangent vectors on $P$.

For each $x\in M$, the fibre ${\ad(P)}_x$ of ${\ad(P)}$ over $x$ has a
Lie algebra
structure (it is actually naturally isomorphic to the Lie algebra of
the isometry group of the fibre $Z_x$). Every element $\a$ of ${\ad(P)}_x$
induces a Killing vector field denoted $\a^*$ on the corresponding
fibre $Z_x$. If $\a$ is represented by $A\in \gg$ in the frame $u\in
P_x$ (i.e. $\a=[u,A]$), and $z\in Z_x$ is represented by $f\in F$ in
the same frame $u$ (i.e. $z=[u,f]$),
then $\a^*_z$ is the image of $A$ by the differential at the identity
of the map $G\to Z_x$, $a\mapsto [u,af]$. By a slight abuse of
notation, we denote this by $\a^*_z=uAf$. It is easy to check that
this is independent of $u$: If we replace $u$ by $ug$, then
$\a=[ug,\ad_{g^{-1}}(A)]$, $z=[ug,g^{-1}f]$, so
$\a^*_z=ug(g^{-1}Ag)(g^{-1}f)=uAf$.

Every section $\a$ of ${\ad(P)}$ induces in this way a vertical
vector field $\a^*$ on $Z$. 
\begin{ede} The vertical vector fields on $Z$ obtained in this way from
sections of ${\ad(P)}$, and the horizontal lifts $X^*$ 
of vector fields $X$ on $M$ are called
{\em standard vertical and horizontal} vector fields on $Z$.
\end{ede}
We recall the classical formulas giving the Lie brackets of standard vertical or
horizontal vector fields on a principal fibration in terms of the
covariant derivative and its curvature form (cf. \cite{kn1}, Ch. 2,
Section 5 or \cite{gau}, Equations (3.9) and (4.4)):
\begin{elem}\label{lb}
If $X,Y$ are vector fields on $M$ and $\a$ is a section of ${\ad(P)}$, then
\beq\label{lb1} [X^*,\a^*]=(\n^{\theta}_X\a)^*,
\eeq
and
\beq\label{lb2} [X^*,Y^*]=[X,Y]^*-\o(X,Y)^*,
\eeq
where $\n^{\theta}$ is the covariant derivative on ${\ad(P)}$ induced by the
connection $\theta$ on $P$ defined in Proposition $\ref{conn}$ and $\o$
is the curvature of $\theta$, viewed as a $2$-form on $M$ with values in
${\ad(P)}$.
\end{elem}

Formula (\ref{lb2}) is equivalent to the fact that  we see that O'Neill's
tensor $A$ associated to the Riemannian submersion $Z\to M$ is given
by $A(X^*,Y^*) = -\frac12\omega (X,Y)^*$ for every vector fields $X, Y$ on $M$
(cf. \cite{besse}, Definition 9.20 and Proposition 9.24).

Using formulas (9.28e) and (9.28c) in \cite{besse} we thus obtain:
\beq\label{e1}
g_Z(R^Z_{X^*,\a^*}Y^*,T^*)=\frac12g_Z(\a^*,((\n^{\theta}_{X}\o )(Y,T))^*),
\eeq
\beq\label{e2}\begin{split} g_Z(R^Z_{X^*,\a^*}Y^*,\b^*)=
  &-\frac14\sum_{a=1}^ng_Z(\a^*,\o
(Y,X_a)^*)g_Z(\b^*,\o (X,X_a)^*) \\&+\frac12g_Z(\b^*,\n^Z_{\a^*}\o
(X,Y)^*).\end{split}
\eeq

\subsection{Curvature constancy}\label{cco}

Let $(Z,g_Z)$ be a Riemannian manifold. For every $z\in Z$ we define
the curvature constancy at $z$ by (see \cite{gray}):
\beq\label{cn}
\V_z:=\{V\in T_zZ\ |\ R^Z_{V,X}Y= g_Z(X,Y)V - g_Z(V,Y)X \
\hbox{for every }X,Y\in T_zZ\}.
\eeq
The function $z\mapsto \dim(\V_z)$ is upper semi-continuous on
$Z$. By replacing $Z$ with the open subset where this function attains 
its minimum, we may assume that $\V$ is a
$k$-dimensional distribution on $Z$, called the
\textit{curvature constancy}.
It is easy to check that $\V$ is totally geodesic
(cf. \cite{gray}).

We will introduce the following Ansatz in order to study the 
curvature constancy condition: Assume that $\V$ is locally 
the vertical distribution of a
Riemannian submersion $\pi: Z\to M$ (equivalently, the metric of $Z$
is bundle-like along $\V$). Since $\V$ is totally geodesic, the fibres of the Riemannian
submersion are locally isometric to the unit sphere $\SM^k$.
All computations below being local, we can assume, by restricting
to a contractible neighbourhood $M'$ of $M$ and taking the universal
cover of $\pi ^{-1}(M')$, that each fibre is globally isometric
to $\SM^k$. Consider the $G$-principal fibre bundle $P$ over $M$
defined in the previous subsection, together with the
connection $\theta$ given by Proposition \ref{conn}.
We set $k+1=:\rr$ so $G=\SO(\rr)$, and
introduce the rank $\rr$ Euclidean vector bundle $E\to M$ associated
to $P$ via the
standard representation of $\SO(\rr)$. Notice that ${\ad(P)}$ is naturally
identified with the bundle $\End^-(E)$ of skew-symmetric endomorphisms
of $E$, and $Z$ is identified with the unit sphere
bundle of $E$.

The curvature constancy condition \eqref{cn} can be expressed in terms
of standard vertical and horizontal vector fields as follows:
\beq\label{cn1}
R^Z_{\a ^*,X^*}Y^*= g_M(X,Y)\a ^*  \quad \hbox{for every }X,Y\in TM,\
\a\in\ad(P).
\eeq
Using \eqref{e1} and \eqref{e2}, this is equivalent to the system
\beq\label{e1p}
(\n^{\theta}_{X}\o )(Y,T)=0,\quad \hbox{for all}\ X,Y,T\in TM, \eeq
\beq\label{cn2}\begin{split}
  g_M(X,Y)g_Z(\a^*,\b^*)=&\frac14\sum_{a=1}^ng_Z(\a^*,\o
(Y,X_a)^*)g_Z(\b^*,\o (X,X_a)^*) \\&-\frac12g_Z(\b^*,\n^Z_{\a^*}\o
(X,Y)^*),\end{split} \eeq 
for all $X,Y\in TM$ and $\a,\b\in\ad(P).$
In order to exploit \eqref{cn2}, we need to express the scalar
product and covariant derivative of standard vertical vector fields
in terms of the corresponding objects on $E$.
\begin{elem}\label{sn}
For every $z\in Z\subset E$ and
$\a,\b,\c\in\ad(P)=\End^-(E)$ in the fibre over $x:=\pi(z)$ we have
\beq\label{sp1}g_Z(\a^*,\b^*)_z=g_E(\a z,\b z).
\eeq
\beq\label{sp2}g_Z(\n^Z_{\a^*}\gamma^*,\b ^*)_z=g_E(\gamma \a z,\b z).
\eeq
\end{elem}
\begin{proof}
Any frame $u$ of $P$ defines an isometry from $(E_x,g_E)$ to the standard
Euclidean space $\RM^\rr$. Once we fix such a frame, $\ad(P)_x$ becomes
the space of skew-symmetric matrices, $Z_x$ is the unit sphere in
$\RM^\rr$, and the vertical vector field $\a ^*$ associated to a skew-symmetric
matrix $\a\in \so(\rr)$ is the Killing
vector field on $\SM^{\rr-1}$ whose value at $z\in \SM^{\rr-1}\subset
\RM^{\rr}$ is $\a z\in T_z\SM^{\rr-1}$. The first formula is now clear.

The Levi-Civita covariant derivative on $\SM^{\rr-1}$ is the projection of the
directional derivative in $\RM^{\rr}$. Moreover, the derivative of the
vector-valued function $f(z)=z$ on $\RM^{\rr}$ obviously satisfies $A.f=A$
for every tangent vector $A\in T\RM^{\rr}$.
We thus get at $z$:
$$
g_Z(\n^Z_{\a^*}\gamma^*,\b
^*)_z=g_E(\a z.\c f,\b z)=g_E(\c \a z,\b z).
$$
\r

Taking Lemma \ref{sn} into account,
\eqref{cn2} is equivalent to
\beq\label{e3}\begin{split}g_M(X,Y)g_E(\a z,\b z)=&\frac14\sum_{a=1}^n
g_E(\a z,\o(Y,X_a) z)
g_E(\b z,\o(X,X_a)z)\\&-\frac12g_E(\b z,\o(X,Y)\a z),\end{split}
\eeq
for all $z\in Z=S(E)$, $\a,\b\in \ad(P)=\End^-(E)$ and $X,Y\in TM$.

Formula \eqref{e3} can be equivalently
stated as follows: 
\beq\label{e4}g_M(X,Y)g_E(v_1,v_2)=\frac 14\sum_{a=1}^n
g_E(v_1,\o(Y,X_a)u)g_E(v_2,\o(X,X_a)u)
-\frac12g_E(v_2,\o(X,Y)v_1),
\eeq
for all $u,v_1,v_2\in E_x$ with $|u|_E ^2=1$ and $v_1,v_2\perp u$ and
for all $X,Y\in 
T_xM$. We introduce the map $\f:\L^2E\to\End^-(TM)$, defined by
$$g_M(\f(u\wedge v)X,Y):=-\frac12g_E(v,\o(X,Y)u),\qquad \forall\
x\in M,\ u,v\in E_x,\ X,Y\in T_xM.$$ 
Formula \eqref{e4} is then equivalent to
\beq\label{e5}\f(u\wedge v)\circ \f(u\wedge w)=\f(v\wedge w)-g_E(v,w)\id,
\eeq
for all $u,v,w\in E_x$ with $|u|_E ^2=1$ and $v,w\perp u$ (where $\id$
denotes the identity of $T_xM$). 

Using the universality property
of the even Clifford algebra (Lemma \ref{up} below), this shows that 
$(E,\f)$ defines an even Clifford structure on $M$. 
We have proved the following:

\begin{ath}\label{3.1}
Assume that the curvature constancy of $Z$ is the vertical distribution
of a Riemannian submersion $(Z^{k+n},g_Z)\to (M^n,g)$. Then $(M,g)$
\begin{itemize}
\item[(a)]
carries a parallel even Clifford structure $(E,\n^E,\f)$ of rank
$\rr=k+1$;
\item[(b)]
the curvature of $E$, viewed as an endomorphism
$\o:\L^2(TM)\to\End^-(E)$, equals minus twice the metric adjoint of
$\f:\L^2E\simeq\End^-(E)\to\End^-(TM)\simeq\L^2(TM)$.
\end{itemize}
Conversely, if $(M,g)$ satisfies these conditions,
then the sphere bundle $Z$ of $E$, together with the Riemannian metric
induced by the connection $\n^E$ on $Z$ defines a Riemannian
submersion onto $(M,g)$ whose vertical distribution belongs to the
curvature constancy.
\end{ath}


\vs

\subsection{The classification}\label{cla}

From Theorem \ref{3.1}, every Riemannian submersion $(Z^{k+n},g_Z)\to
(M^n,g)$ whose vertical distribution belongs to the
curvature constancy defines a parallel even Clifford
structure $(E,\n^E,h,\f)$ of rank $\rr:=k+1$ on $M$, such that the
curvature $\o$ of $\n^E$, viewed as an endomorphism
$\o:\L^2(TM)\to\End^-(E)$, equals minus twice the metric adjoint of the
Clifford morphism
$\f:\L^2E\to\End^-(TM)$. In the notation of Proposition \ref{p37}, this
amounts to say that
\beq\label{of} \o_{ij}=2J_{ij},\qquad\forall \ 1\le i\ne j\le \rr.
\eeq
Conversely, if $(E,\n^E,h,\f)$ is a parallel even Clifford
structure of rank $\rr$ on $M$ satisfying \eqref{of}, $E$ carries a
Riemannian metric defined by the metric on $M$, that of $E$, and the
splitting of the tangent bundle of $E$ given by the connection
$\n^E$ and by Theorem \ref{3.1}, the restriction to the unit sphere
bundle $Z$ of the projection $E\to M$ is a Riemannian submersion
whose vertical distribution belongs to the curvature constancy.

We will now examine under which circumstances a simply connected
complete Riemannian manifold
\begin{enumerate}
\item[(i)] carries a parallel even Clifford structure $(E,\n^E,h,\f)$.
\item[(ii)] $(E,\n^E,h,\f)$ satisfies \eqref{of}.
\end{enumerate}
Notice that for every $3\le n\ne 4$, condition (ii) together with
\eqref{h} and \eqref{0r2} implies that the scalar curvature of $M$ is
\beq\label{scal}\scal=2n(n/4+2\rr-4).
\eeq

\noindent$\bullet$ {$\rr=2$.} In this case $M$ is K\"ahler (see
Example \ref{4.1}) and $E$ is simply a rank 2
Euclidean vector bundle endowed with a metric connection $\n^E$ whose
curvature is minus twice the K\"ahler form of $M$. By the Chern-Weil theory,
this is equivalent to the cohomology class of the K\"ahler form being
half-integer, so up to rescaling
$M$ is a Hodge manifold. It is well known that the circle bundle $Z$
of $E$ carries a Sasakian structure for the corresponding rescaling
of the metric on $M$.

\noindent$\bullet$ {$\rr=3$.} By Example \ref{4.2}, condition (i) is
equivalent to $M$ being quaternion-K\"ahler (recall that this is an
empty condition for $n=4$) and $E$ is either $\L^2_+M$ for $n=4$ or
the 3-dimensional sub-bundle of $\L^2M$ defining the
quaternion-K\"ahler structure for $n>4$. Condition (ii) is
equivalent to $M$ being anti-self-dual and Einstein with scalar
curvature equal to 24 (see \cite{besse} p.51 and \eqref{scal} above)
for $n=4$, and quaternion-K\"ahler with positive scalar curvature
equal to $8q(q+2)$ for $n=4q>4$. The Riemannian manifold $(Z,g_Z)$
is the twistor space of $M$ in the sense of Salamon \cite{sal}.

\noindent$\bullet$ {$\rr=4$.} Proposition \ref{p37} (i) shows that $M$
is the Riemannian product of two quaternion-K\"ahler manifolds $M^+$
and $M^-$ of dimension $4q^+$ and $4q^-$ respectively (notice that one
of $q^+$ or $q^-$ might vanish). Recall that the rank 4 even Clifford structure
$E$ on $M$ induces in a natural way rank 3 even Clifford structures
$\L^2_\pm E$ on $M^\mp$. A local
orthonormal basis $e_i$, $1\le i\le 4$ of $E$ induces local
orthonormal bases $\t e ^\pm_i$, $1\le i\le 3$ of $\L^2_\pm E$ by
\eqref{te}. Taking \eqref{tj} and \eqref{o} into account, Equation
\eqref{of} becomes
\beq\label{of1} \o ^\pm_{ij}=4J^\pm_{ij},\qquad\hbox{for all} \ 1\le
i\ne j\le 3. 
\eeq
Like in the previous case, this means that $M^\pm$ is a
quaternion-K\"ahler manifold with scalar curvature $16q^\pm(q^\pm+2)$,
where now we use the usual convention that in dimension 4
quaternion-K\"ahler means anti-self-dual and Einstein.

In order to
describe the Riemannian manifold $(Z,g_Z)$, we need to understand in
more detail the construction of the even Clifford structure of rank 4
on a product of quaternion-K\"ahler manifolds $M=M^+\times
M^-$. The orthonormal frame bundle of $M$ admits a reduction to a
principal bundle $P$ with structure group $G:=\Sp(q^+)\.\Sp(1)\times
\Sp(q^-)\.\Sp(1)\subset \SO(4q^++4q^-).$ The representation of the
universal cover of $G$ on $\RM^4\simeq \HM$ defined by
$$(A,a,B,b)(v)=avb^{-1}\qquad \forall\ (A,a,B,b)\in\Sp(q^+)\times\Sp(1)\times
\Sp(q^-)\times\Sp(1),\ \forall v\in \HM$$
induces a projective representation $\rho:G\to\PSO(4)$, which in turn
determines the (locally defined) bundle $E$ and the (globally defined)
manifold $Z:=P\times_\rho\RM\PM^3$. A Riemannian manifold obtained in this
way is called {\em quaternion-Sasakian}. By definition, a quaternion-Sasakian
manifold fibres over a product of quaternion-K\"ahler manifolds $M=M^+\times
M^-$, with fiber $\RM\PM^3$. Notice that 3-Sasakian
manifolds are special cases of  quaternion-Sasakian manifolds, when
one of the factors $M^+$ or $M^-$ is reduced to a point.

We now examine the remaining cases in Table 2.

\noindent$\bullet$ {$\rr\ge 5$ and $n=8$.} Taking \eqref{2r} into
account, \eqref{of} is equivalent to the fact that the restriction of
the curvature endomorphism $R$ of $M$ to the Lie sub-algebra
$\f(\L^2E)\subset\L^2M$ equals $4\id$. Moreover, we have $\kk=2$ in
Equation \eqref{h}, so \eqref{0r1} shows that $M$ is Einstein with
scalar curvature $2n(n/4+2\rr-4)$.

If $\rr=8$, this means that $R$ is constant, equal to 4 on $\L^2M$, so
$M$ is the round sphere $\SM^8(1/2)$ of radius 1/2.

The case $\rr=7$ does not occur, since a manifold with holonomy
$\Spin(7)$ is Ricci-flat, contradicting Proposition \ref{p37} (iii).

If $\rr=6$, $M$ is K\"ahler and $\f(\L^2E)$ is just the sub-bundle
$\L^{(1,1)}_0M$ of primitive forms of type $(1,1)$, corresponding to
the isomorphism $\spin(6)\simeq \su(4)$. By the above $M$ has Einstein
constant equal to 20. This shows that the curvature endomorphism of
$M$ is equal to 4 on $\L^{(1,1)}_0M$, is equal to 20 on the line
generated by the K\"ahler form (since the image of the K\"ahler form
is the Ricci form), and vanishes on $\L^{(2,0)+(0,2)}M$ (like on every
K\"ahler manifold), so $M$ is isometric to the complex projective
space $\CM\PM^4=\SU(5)/{\rm S}(\U(1)\.\U(4))$ endowed with the
Fubini-Study metric with scalar curvature 160.

If $\rr=5$, $M$ is quaternion-K\"ahler, and by a slight abuse of
notation we can write $\L^2M=\sp(1)\oplus\sp(2)\oplus\pp$. Like
before, the curvature endomorphism $R$ of $M$ equals 4 on
$\sp(2)=\f(\L^2E)$. Moreover, on every quaternion-K\"ahler manifold
with Einstein constant 16, $R$ equals 4 on $\sp(1)$ and vanishes on
$\pp$. Thus $M$ is isometric to $\HM\PM^2=\Sp(3)/\Sp(1)\times\Sp(2)$.

\noindent$\bullet$ {$\rr\ge 5$ and $n>8$.} This case concerns the
symmetric spaces $M$ in the last seven rows of Table 2. For each of
these spaces condition (ii) is automatically satisfied (by Proposition
\ref{p37}) for the specific normalization of the metric for which
$\kk=2$ in Equation \eqref{h}, which by \eqref{0r1} is equivalent to
the scalar curvature being equal to $2n(n/4+2\rr-4)$.

Summarizing, we have proved the following
\begin{ath}\label{fc}
There exists a Riemannian submersion
from a complete simply connected Riemannian manifold $(Z^{k+n},g_Z)$
to a complete simply connected Riemannian manifold
$(M^n,g)$ whose vertical distribution belongs to the
curvature constancy if and only if $(Z,M)$
appears to the following list:
\end{ath}
\vs
\begin{tabular}{|c|l|l|r|r|}\hline
  $Z$  &  $M$   &  Fibre  & $\dim(M)$    &  $\scal(M)$  \\

 \hline\hline
Sasakian       &   Hodge   &   $\SM^1$  &   $2m,\ m\ge1$ &       \\
\hline
Twistor space   $Z$  &   quaternion-K\"ahler (QK)& $\SM^2$ &$4q,\ q\ge1$  &  $
8q(q+2) $ \\
\hline
Quaternion-Sasakian    &  product of two QK &   $\RM\PM^3$ &
$4(q^++q^-),$ & $16q^+(q^++2) $ \\
           & manifolds &    &  $q^++q^-\ge 1 $ &$ +16q^-(q^-+2)$   \\
\hline
$\frac{\Sp(q^++1)\times\Sp(q^-+1) }{\Sp(q^+)\times\Sp(q^-)\times \Sp(1)} $  & $\HM\PM^{q^+}\times \HM\PM^{q^-}$  &   $\SM^3$ &
$4(q^++q^-),$ & $16q^+(q^++2) $ \\
           &     &    &  $q^++q^-\ge 1 $  &$ +16q^-(q^-+2)$   \\
\hline\hline
$\frac{\Sp(k+2)}{\Sp(k)\times\Spin(4)}$ & $\Sp(k+2)/\Sp(k)\times\Sp(2)$ & $\SM^4$  &
$8k,\ k\ge1$  &  $32k(k+3)$ \\
\hline
$\frac{\SU(k+4)}{{\rm S}(\U(k)\times(\Sp(2)\.U(1)))}$ & $\SU(k+4)/{\rm
  S}(\U(k)\times\U(4))$ &$\RM\PM^5$
& $8k,\ k\ge1$ &   $32k(k+4) $\\
\hline
$\frac{\SO(k+8)}{\SO(k)\times\Spin(7)}$ &
$\SO(k+8)/\SO(k)\times\SO(8)$ &   $\RM\PM^7$ &
$8k,\ k$ odd $\ge3$  &  $32k(k+6) $ \\
\hline
$\frac{\Spin(k+8)}{\SO(k)\times\Spin(7)}$ &
$\SO(k+8)/\SO(k)\times\SO(8)$ &   $\SM^7$ &
$8k,$ $k=1$ or  &  $32k(k+6) $ \\
 &
 &   &
 $k$ even &  \\
\hline\hline
$F_4/\Spin(8)$   & $F_4/\Spin(9)$ &  $\SM^8$ &  16 & $ 2^6\.3^2 $ \\
\hline
$ E_6/\Spin(9)\.\U(1)$ & $ E_6/\Spin(10)\.\U(1)$ &  $\SM^9$   &  32&
$2^9\.3 $\\
\hline
$ E_7/\Spin(11)\.\SU(2)$  & $ E_7/\Spin(12)\.\SU(2)$ &  $\SM^{11}$   &
64&  $2^9\.3^2 $ \\
\hline
$E_8/\Spin(15)$   & $E_8/\Spin^+(16)$  &  $\SM^{15}$    &  128 &
$2^{10}\.3\.5 $ \\
\hline
\end{tabular}
\vskip.2cm

\begin{center} Table 3. Riemannian submersions with curvature
  constancy.\footnote{We adopt in this table the
  usual convention for quaternion-K\"ahler manifolds in dimension 4 as
  being anti-self-dual and Einstein.}
\end{center}

\vs

In particular, the above table shows that all Hopf fibrations provide
examples of manifolds with curvature constancy.

We end up this section with a short list of interesting problems 
related to Clifford structures and perspectives of 
possible further research. These are just a few examples of the
numerous questions raised by our work.
\begin{itemize}
\item The notion of curvature constancy has a hyperbolic counterpart which 
leads to the notion of Lorentzian Clifford structures. This problem 
can be studied with methods similar to those above and could 
provide a new framework for theoretical physicists.
\item Many notions and results from almost Hermitian geometry
can be generalized to Clifford structures. One can for instance
introduce the
minimal connection of an (even) Clifford structure, and obtain 
a Gray-Hervella-type classification of Clifford structures.
\item
One can also address the question of the existence 
of global almost complex structures compatible with a parallel (even)
Clifford structure, generalizing corresponding results by D. V. Alekseevsky, 
S. Marchiafava and M. Pontecorvo \cite{AMP} obtained for
quaternion-K\"ahler manifolds.
\end{itemize}

\vs\vs

\section{Appendix. The universality property of the even Clifford algebra}

For the reader's convenience we provide here the proof of the
universality property for even Clifford algebras which was needed in
the proof of Theorem \ref{fc}.

\begin{elem}\label{up}
Let $(V,h)$ be a Euclidean vector space and let $\A$ be any real
algebra with unit. 
We make the usual identification of $\L^2V$ with a subspace
of $\Cl^0(V,h)$.
Then a linear map $\f:\L^2V\to\A$ extends to an algebra morphism
$\f:\Cl^0(V,h)\to\A$ if and only if it
satisfies
\beq\label{e51}\f(u\wedge v)\circ \f(u\wedge w)=\f(v\wedge w)-h(v,w)1_{\A}
\eeq
for all $u,v,w\in V$ with $|u|_V^2=1$ and $v,w\perp u$. 
\end{elem}
\bp
The ``only if" part is obvious. Assume, conversely, that 
\eqref{e51} holds and let $u,v,w\in V$ be arbitrary vectors. 
We apply \eqref{e51} to the triple $\t u:=u/|u|$,
$\t v:=v-h(u,v)u/|u|^2$ and $\t w:=w-h(u,w)u/|u|^2$ and obtain
\beq\label{e511}\begin{split}\f(u\wedge v)\circ \f(u\wedge
  w)=&|u|^2\f(v\wedge w) 
-h(u,v)\f(u\wedge w)- 
h(u,w)\f(v\wedge u)\\&-(|u|^2h(v,w)-h(u,v)h(u,w))1_{\A}. \end{split}
\eeq
By defining
$\s:V\otimes V\to\A, 
u\otimes v \mapsto \s_{uv}:=\f(u\wedge v)-h(u,v)1_{\A}$, \eqref{e511}
becomes equivalent to
\beaa
\label{s1}\s_{uv}+\s_{vu}&=&-2h(u,v)1_{\A},\\
\label{s2}\s_{vu}\circ\s_{uw}&=&-h(u,u)\s_{vw}
\eeaa
for all $u,v,w\in V$.
Let $T(V)$ denote the tensor algebra of $V$ and
$$T^0(V):=\bigoplus_{k\ge0}V^{\otimes 2k}.$$
By definition, $\Cl^0(V,h)=T^0(V)/\I$, where $\I$ is the intersection
with $T^0(V)$ of the two-sided ideal of $T(V)$ generated by elements
of the form $u\otimes u+h(u,u)$. The map $\s$ clearly induces a unique
algebra morphism $\s^*:T^0(V)\to\A$ such that $\s^*=\s$ on
$V\otimes V$. We claim that $\I\subset\Ker(\s^*)$.
Now, every element of $\I$ is a linear combination of elements of the
form $A=a\otimes(u\otimes u+h(u,u))\otimes b$ or $B=a\otimes
v\otimes(u\otimes u+h(u,u))\otimes w\otimes b$, with $a,b\in T^0(V)$
and $u,v,w\in V$. From \eqref{s1} we have
$$\s^*(A)=\s^*(a)\circ(\s_{uu}+h(u,u)1_{\A})\circ\s^*(b)=0,$$
and \eqref{s2} yields
\bea\s^*(B)&=&\s^*(a)\circ\s^*(v\otimes u\otimes u\otimes
w+h(u,u)v\otimes w)\circ\s^*(b)\\
&=&\s^*(a)\circ(\s_{vu}\circ\s_{uw}+h(u,u)\s_{vw})\circ\s^*(b)=0.\eea
Consequently $\s^*$ descends to an algebra morphism $\Cl^0(V,h)\to\A$,
whose restriction to $\L^2V$ is just $\f$.
\r


\labelsep .5cm

\end{document}